\documentclass[12pt]{article}  

\usepackage{amsmath}  
\usepackage{amsfonts}  

\usepackage{amssymb}  
\usepackage[all]{xy}  

\usepackage[mathscr]{euscript}  
\renewcommand\L{\Lambda}  
\newcommand\F{\mathcal F}  
\newcommand\NN{\mathcal N}  
\newcommand\G{\mathbf G}  
\renewcommand\H{\mathcal H}  
\newcommand\J{\mathbf F}  
\newcommand\f{\mathbf f}  
\newcommand\g{\mathbf g}  
  
\newcommand\D{\mathcal D}  
\newcommand\LL{\mathcal L}  
\renewcommand\l{\lambda}  
\newcommand\C{\mathbb{C}}  
\newcommand\R{\mathbb{R}}  
\newcommand\Z{\mathbb{Z}}  
\newcommand\N{\mathbb{N}}  
\renewcommand\P{\mathbb{P}}  
\newcommand\Q{\mathbb{Q}}  
\newcommand\QQ{\mathcal Q}  
\newcommand\ev{\operatorname{ev}}  
  
\newcommand{\m}  {\mathbf m}  

\newcommand\p{\partial}  
\newcommand\eps{\varepsilon}  
  
\newcommand\Aut{\operatorname{Aut}}  
  
\renewcommand\a{\alpha}  
\renewcommand\b{\beta}  
\renewcommand\c{\gamma}  
\newcommand\Res{\operatorname{Res}}  
\newcommand\w{\wedge}  
\newcommand\Lie{\operatorname{Lie}}

\renewcommand\o{\overline}  
\newcommand\V{\mathcal V}

\begin{document}  

\title{\let\thefootnote\relax\footnotetext{This work was supported by IBS-R003-D1.}  A Mirror Theorem for $T $-Equivariant Blowups}  
\author{Jeff Brown}  
\date{}  

\maketitle
\begin{abstract}  

Let $E $ be a toric fibration arising from symplectic reduction of  a  direct sum of  line bundles over (almost-) K\"ahler base $B $.  Then each torus-fixed point of the toric manifold fiber defines a section of the fibration.  Let $L_a $ be convex line bundles over $B $, $A_a $ smooth divisors of $B $ arising as the zero loci of generic sections of $L_a $, and  $\a:B\to E $  a particular fixed-point section of $E $.  Further assume the $\{A_a\}    $ to be  mutually disjoint.

We compute genus-0 Gromov--Witten invariants of 
the blowup of $E $ along $\a(\coprod_a A_a)  $ in terms of genus-0 Gromov--Witten invariants of $B $ and of $\{ A_a\} $, the matrix used for the symplectic reduction description of the fiber of the toric fibration $E\to B $, and the restriction maps $i_{A_a}  ^*:H^*(B)   \to H^*(A_a)  $.

Let $\sum_DQ^DJ_{e(\ \cdot\ )   , L_a}  ^D $ be the $Euler(L_a)  $-twisted J-function of the base, where the summation index runs over all curve classes    in the image of the 
  push-forward ${i_{A_a}}  _*:H_2(A_a, \Z)   \to H_2(B, \Z)  $.  
The Quantum Lefschetz Theorem \cite{Coates-Givental}   relates it to the genus-0 Gromov--Witten invariants of the complete intersection $A_a $.
 The  geometry of the present blowup gives (indirectly, at least)    a Gromov--Witten  theoretic interpretation to $J_{e(\ \cdot\ )   , L_a}  ^D $ {\em for all  curve classes $D $ in $H_2(B)  $}.
\end{abstract}  

\section{Formulations}  
                        {\bf 1.1. Genus-0 Gromov--Witten invariants.}   Associated to an
                        almost-K\"ahler manifold $M $ are the moduli spaces $M_{0, n, D} $ of (equivalence classes of)    degree-$D $ stable
                      maps into $M $ of genus-$0 $ (possibly nodal)    compact connected
                        holomorphic curves with $n $ marked points.   Two such stable maps $(f;C, x_1, \dots, x_n)  $ and $(f';C', x'_1, \dots, x'_n)  $ are {\em equivalent}   if there is a holomorphic automorphism $\phi:C\to C' $ mapping marked points to marked points and preserving the ordering, such that  $f=f'\phi $.   For a stable map $(f;C, x_1, \dots, x_n)  $, the degree-$D $ condition reads $f_*([C])   =D $.  Define the {\em genus-0 descendant
                        potential of $M $}   to be the formal series  
                        \[ \F_M(t)   =\sum_{n=0}  ^{ \infty}  \sum_{D\in MC}  \frac{Q^D}  {n!}  \int_{[M_{0, n, D}  ]}  \prod_{a=1}  ^n\sum_{k=0}  ^{\infty}  \ev_a^*(t_k)   \psi_a^k, 
                      \]
                        where $[M_{0, n, D}  ]$ denotes  the virtual  fundamental class of $M_{0, n, D} $, $MC $---the Mori cone of $M $, that   is the semigroup in
                       $H_2(M, \Z)  $ generated by classes representable by  compact holomorphic
                        curves, $Q^D $---the element in the {\em Novikov ring}  
                        (that   is a power-series  completion of the semigroup algebra of the Mori cone)   
                        representing the degree $D \in MC $, $\psi_a $ ---
                        the 1st Chern class of the universal cotangent line bundle over $M_{0, n, D} $ formed by
                        the cotangent
                        lines along
                         the stable
                      maps at the $a $-th marked point, $\ev_a $ ---
                        the map $M_{0, n, D}  \to M $
                        that  evaluates the  stable maps at the 
                       $a $-th marked point, $t_k \in H^*(M, \QQ)   , \
                        k=0, 1, 2, \dots, $---arbitrary cohomology classes of $M $
                        with coefficients in a suitable {\em ground ring}  $\QQ $ (for the moment
                        let it be the rational Novikov ring $\Q [[MC]]$).
                        
                        {\bf 1.2. Toric fibrations.}   Let $\m=(m_{ij}  |i=1, \dots, K;
                        j=1, \dots, N)  $ be an integer matrix, and consider the action of $T^K $
                        on the Hermitian space $\C^N $ that, for each $j=1, \dots, N $, multiplies the coordinate $z_j $ by \linebreak
                       $\exp(\sum_{i=1}  ^Km_{ij}  \sqrt{-1}  \theta_i)  $.
Let $\mu: \C^N\rightarrow \R^N $ be the map given by \[ (z_1, z_2, \dots, z_N)   \to (|z_1|^2, |z_2|^2, \dots, |z_N|^2)   .\] If the {\em moment map}  $\m\circ\mu $ has a regular value $\omega\in\R^K $, then $(\m\circ\mu)   ^{-1}  (\omega)   /T^K $ is a symplectic manifold.  This construction is called {\em  symplectic reduction}.  The space $(\m\circ\mu)   ^{-1}  (\omega)   /T^K $ is also denoted by $\C^N//_{\omega}   T^K $, and isequipped with a canonical symplectic form, call it $\omega $, induced by the standard symplectic form on $\C^N $. Given complex line bundles $L_1, \dots, L_N $ form the vector bundle $\oplus      L_j\rightarrow
B $.   Let us assume that  $T:=T^N $ is the structure group of the vector
   bundle $\oplus   L_j $ as follows. All complex line bundles over $B $ may be assumed to have the unitary circle $S^1 $ as structure group, as they are induced by pullback from the tautological line bundle over $\C P^{\infty} $.  Since the fiberwise moment map is $T $-invariant it follows that   the fiberwise symplectic reduction 
of $\oplus   L_j $ is well-defined giving the toric fibration $E\rightarrow B $.
The $i $th coordinate $\theta_i $ on the torus $T^K $ defines a circle bundle over $E $ for which the expression $\sqrt{-1}  d\theta_i $ defines connection 1-forms in the bundle.  Denote by $-P_i $ the first Chern class of the $i $th circle bundle over $E $, and $-p_i $ its restriction  to a fiber.  The $p_1, \dots, p_K $ classes are of Hodge $(1, 1)  $-type  by the Fubini--Study construction, though they need not be K\"ahler classes\footnote{Section 1.6 gives a description of a toric manifold for which the class $p_3 $ is non-K\"ahler.}.   Let $\c $ be any $T $-fixed point of $(\m\circ \mu)   ^{-1}  (\omega)/T^K  $, and $p\in(\m\circ\mu)^{-1}  (\omega)  $
representing the $T^K $-equivalence class $\c $.  The orbits $T^Kp $ and $Tp $ are then identical.  It follows that  there is some coordinate subspace $\C^K\subset \C^N $ with coordinates $z_{j_1}  , \dots, z_{j_K} $, containing $p $, such that   none of the coordinates $z_{j_1}  (p)   , \dots, z_{j_K}  (p)  $ vanishes.   It will be convenient to think of the $T $-fixed strata $\c $ of $E $ in terms of the corresponding indices $j_1, \dots, j_K $.   Define $-\L_j:=c_1^T(L_j)  $  for $j=1, \dots, N $.  For each $j=1, \dots, N $, the restriction of $U_j=\sum_{i=1}  ^Km_{ij}  P_i-\L_j $ to a fiber is Poincar\'e dual to the $j $th coordinate divisor $((\m\circ \mu)   ^{-1}  (\omega)   \cap \{z_j=0\}    )   /T^K $.  The expressions for the pullbacks $P^{\c}  _i $ in terms of $\L_j $ may be summarized by the equations $\c^*U_{j_1}  =\cdots =\c^*U_{j_K}  =0 $.
All bundles introduced thus far are $T $-equivariant, so their Chern classes may be assumed to take values in the $T $-equivariant cohomology group $H_T^*(E) $ with coefficient ring $H^*(BT, \Q)=\Q[\l_1, \dots, \l_N]$.

{\bf 1.3. The cone $\LL_{E_{\a(A) }} $.}   Associated to the genus-0 Gromov--Witten theory of $M $ is a {\em Lagrangian cone}  $\LL_M $ in a symplectic loop space $(\H, \Omega)  $ \cite{Coates-Givental}.
The space $\H=\H_+\oplus \H_- $ is a module over the  ground ring $\QQ $. Pending further completions, $\H $
consists of Laurent series  in $1/z $ with coefficients in $H:=H^*(M, \QQ)  $, completed so that $\H_+ $ consisting of elements of $H [z]$ {\em at each order  in Novikov's variables}  , and $\H_-:=z^{-1}  H [[z^{-1}  ]]$.  Identify each $q(z)   =\sum_{k=0}  q_kz^k\in \H_+ $ with the domain variables $t_0, t_1, t_2, \dots $ of $\F_M $
by the {\em dilaton shift}   convention $q_k=t_k-\delta_{k, 1}  , k=0, \dots, \infty $.  Take the ring of coefficients for Novikov's variables to be the (super-commutative)    power series  ring (with coefficients in the field of fractions $\Q(\l):= \Q(\l_1, \dots, \l_N) $, in all of our applications)    in formal
 coordinates along $H^*(M, \Q)  $, and require the  variables $t_0, t_1, \dots $ to vanish
 when Novikov's variables and the newly introduced formal coordinates are all set to zero.  This gives a Novikov ring $\QQ $ that   is consistent with the formula for $I_{E_{\a(A) }} $ in our Main Theorem.
Consider the symplectic manifold $T^*\H_+ $ with standard symplectic form $\sum_{k=0}  ^{\infty}  dp_k\wedge dq_k $.  It is symplectomorphic to $\H $ with symplectic form\[ \Omega(\f, \g)   :=\frac{1}  {2}  \Res_{z=0}  (f(-z)   , g(z)  ) _M, \]
  where $(\cdot, \cdot)_M $ is the Poincar\'e pairing.

Let us implement this symplectomorphism via the map
   
\[(q, p)   \mapsto \sum_{k=0}  ^{\infty}  q_kz^k+ \sum_{k=0}  ^{\infty}  p_k(-z)   ^{-k-1}  .\]

Consider the  graph of the differential of $\F_M(t)  $, which is a Lagrangian submanifold  in $T^*\H_+ $.  From it we arrive at
\[ \LL:=\{(q, p)   | p=d_t\F_M(t)   \}    \] by rigid translation in the direction of the dilaton shift.  Thus $\LL $ is also a Lagrangian submanifold.  Henceforth we consider $\LL $ as a  submanifold of $\H $.  The work of Coates--Givental \cite{Coates-Givental}  establishes that  $\LL $ is a (Lagrangian)    cone with vertex at the $(q, p)  $-coordinate origin; that  is, 
\[ T_{\f}  \LL \cap \LL=zT_{\f}  \LL\ \ \forall \f \in \LL\setminus (0, 0).\]
In particular, each tangent space is  preserved by multiplication by $z $.
 The fact that $\LL $ contains the $(q, p)  $-coordinate origin is a special case 
 of Getzler's \cite{Getzler}, Coates--Givental's \cite{Coates-Givental}    solution of Eguchi--Xiong's, Dubrovin's $(3g-2)  $-jet conjecture. 
 The Lagrangian cone $\LL $ in the 
torus-equivariant genus-0 Gromov--Witten theory of $E_{\a(A) } $ lies in the symplectic
loop space $(\H, \Omega)  $.

A point in the cone can be written as
\begin{align*}   \J (-z, t)    =&\\ -1z+&t(z)+\sum_{n}   \sum_{ D, d, {\tilde d}}  \frac{Q^Dq^d{\tilde q}  ^{\tilde d}}  {n!}  (\ev_1)   _*
\left[ \frac{1}  {-z-\psi_1}  
\prod_{i=2}  ^{n+1}   (\ev_i^*t)   (\psi_i)   \right], \end{align*}  where $(\ev_1)   _* $ denotes the virtual  push-forward by the evaluation
map $\ev_1: (E_{\a(A) }  )   _{0, n+1, \D}  \to E_{\a(A) } $, and $t(z)   =\sum_{k=0}  ^{\infty}  t_kz^k $
is an element of $\H_+ $ with arbitrary coefficients $t_k\in H $.
Define the J-function to be  the restriction of $\J(-z, t)  $ to values $t_0\in H $ and to $t_k=0 $ for all $k>0 $.  For each $\f\in \LL $ there is a unique $t(\f)    \in H $ such that 
\[ zT_{\f}  \LL \cap \{-z+z\H_-\}    =J(-z, t(\f)  )   .\]
 The property of the set of all  tangent spaces of $\LL $ to be in 1-1 correspondence with the set $H $, which is a finite-dimensional $\QQ $-module, is called {\em overruled}.
Let $\{\phi_{\mu}  \}   $ be a basis of $H^*(M)  $ and $\{\phi^{\mu}  \}   $ the Poincar\'e-dual basis.   For each $t\in H $ and for each open set $U\ni t $, the J-function generates a module
over the algebra $\Gamma_U(\oplus _{r=0}  ^{\infty}  ( \otimes^rT_tH)\otimes\Q (z^{-1}  )  )  $ of differential operators as follows, 
\[ z\p_az\p_b J(z, t)   =z\p_{a\bullet_t b}  J(z, t)   , \] where \[a\bullet_t b:= \sum_{n=0}  ^{\infty}  \sum_{D\in MC}  \frac{Q^D}  {n!}  \phi_{\mu}  \int_{[M_{0, n+3, D}  ]}  \ev_1^*\phi^{\mu}  \ev_2^*a  \ev_3^*b \prod_{i=4}  ^{n+3}  \ev_i^*t\]is the unital, associative, (super-)    commutative quantum cup product.   Additionally, the J-function satisfies the {\em string}   and {\em divisor}  equations
\[ z\p_{1_M}  J(z, t)   =J(z, t)   \ \text{and}  \ z\p_{\rho}  J(z, t)   =\sum_D(\rho+\rho(D)z  )J_M^D(z, t) \forall\rho\in H^2(M, \Q)   , \]
respectively.

{\bf 1.4. Twisted Lagrangian cones.}   The forgetful maps $ft_{n+1}  :M_{0, n+1, D}  \to M_{0, n, D} $ induce the K-theoretic  push-forward maps $(ft_{n+1}  )   _*:K(M_{0, n+1, D}  )   \to K(M_{0, n, D}  )  $.   
 Let $\V $ be a complex vector bundle over $M $. The evaluation  maps $\ev_{n+1}  :M_{0, n+1, D}  \to M $ induce the (virtual-)     bundles $\ev^*_{n+1}  \V $, in terms of which the (virtual-)    virtual  bundles
\[ \V_{0, n, D}  :=(ft_{n+1}  )   _*\ev^*_{n+1}  \V\in K(M_{0, n, D}  )   \]are defined.  The fiber of $ \V_{0, n, D} $ over a stable map $(f;\Sigma, p)  $ is \[ H^0(\Sigma, f^*\V)   \ominus
H^1(\Sigma, f^*\V)   .\]
  Given a characteristic class $\bf c(\cdot) $, define the twisted Poincare pairing \linebreak $(a, b)   _{{\bf c(\cdot)}  , \V}  :=(a, {\bf c}  (\V)   b)_M  $.

   A point in the $({\bf c}  (\cdot)   , \V)  $-twisted cone can be written as
\begin{align*}   \J_{{\bf c}  (\ \cdot\ )   , \V}  
(-z, t)    =&\\ -1z+&t(z)+\sum_{n, D}  \frac{Q^D}  {n!}  (\ev_1)   _*
\left[{\bf c}  (\V_{0, n, D}  )   \frac{1}  {-z-\psi_1}  
\prod_{i=2}  ^{n+1}  (\ev_i^*t)   (\psi_i)\right].\end{align*}  The overruled Lagrangian cone $\LL_{{\bf c}  (\cdot)   , \V} $ in the 
 $({\bf c}  (\cdot )   , \V)  $-twisted genus-0 Gromov--Witten theory of $M $ lies in the  symplectic loop space $(\H^{{\bf c}  (\cdot)   , \V}  , \Omega_{{\bf c}  (\cdot\ )   , \V}  )  $, 
where $t(z)   =\sum_{k=0}  ^{\infty}  t_kz^k $
is an element of $\H^{{\bf c}  (\cdot)   , \V}  _+:=(\H_M)_+ $ with arbitrary coefficients $t_k\in H^{{\bf c}  (\cdot)   , \V}  :=H_M $.
The examples we will consider are:

{\bf Example 1.4.1.}  ${\bf c}  (\cdot)   =Euler(\cdot)  $, and $\V $ is a {\em convex}   line bundle; i.e., $H^1 (\Sigma, f^*\V)   =0 $; or equivalently, $f^*c_1(\V)(\Sigma)\geq -1 $  for all genus-0 stable maps $(f;\Sigma, p(\Sigma)  )  $ to $M $.

{\bf Example 1.4.2.}  ${\bf c}  (\cdot)   =Euler_T^{-1}  (\cdot)  $, and $\V $ is a complex vector bundle with a hamiltonian torus action that decomposes $\V $ into a direct sum of complex line bundles, each of which carries a non-trivial $T $-action. 

{\bf 1.5. Torus action on the blowup.}

The blowup ${\tilde \pi}  :E_{\a(A) }  \to E $ along the divisor $\a(A)  $ of the $T $-fixed section
 $\a\subset  E $ may be described as the result of surgery on
 $E $, as we now recall.  

Define a map from
 a tubular neighborhood of the ${\cal O}  (-1)  $ bundle over the projective bundle
 ${\mathbb P}  (N_{\a(A) }  E)  $ over $\a(A)  $ to a tubular neighborhood of 
 $N_{\a(A) }  E $ over $\a(A)  $ as follows.  Fiberwise, it is described by the standard
blowup map 
\[ \{(v, [l])   \in {\C}  ^n\times     {\C P}  ^{n-1}  : \text{ $l $ a line in
 ${\mathbb C}  ^n $ through the origin; $v\in l $}  \}    \rightarrow {\mathbb
C}  ^n, \]
\[ (v, [l])   \mapsto v.\]
This map identifies $T $-equivariantly the complements of the 0-sections of the total spaces of the preceding two vector bundles.  Remove a tubular neighborhood of $N_{\a(A) }  E $ from $E $, and replace it
by a tubular neighborhood of the ${\cal O}  (-1)  $ bundle over
 ${\mathbb P}  (N_{\a(A) }  E)  $.  The $T $-action on $E_{\a(A) } $  induces a $T $-action on the moduli spaces of stable maps to $E_{\a(A) } $, which in turn induces a $T $-action on the universal cotangent line bundles at each of the marked points. We consider $L $, giving rise to $A\subset  B $ (as the zero locus of a generic section), as $T $-equivariant with the trivial $T $-action.  For a given $T $-fixed stratum $\eps $ of $E_{\a(A)} $ and a line bundle $\V \to B $ with a fiberwise $T $-action, we refer to the restriction $\eps^*Euler_T(\V)   \in H^2_T(\eps, \Q)     \simeq  H^2(BT, \Q)   \oplus H^2(\eps, \Q)  $ as the {\em $T $-weight}   of $\V $ at $\eps $.

{\bf Example 1.5.}    If $\eps $ is a $T $-fixed stratum in the complement of the exceptional divisor, then take $M=B\simeq\eps $ and $\V=N_{\eps}  E_{\a(A)} $ in Example 1.4.2. 
If $\eps $ is a $T $-fixed stratum in the exceptional divisor, then  take $M=A\simeq\eps$ in Example 1.4.2 and $\V $ to be the subbundle $\NN^{\eps} $ of $N_{\eps}  E_{\a(A)} $ where $T $ acts non-trivially.

{\bf  1.6. Motivation.}   Let $X $ be a compact symplectic toric manifold and let $T\subset  (\C^*)   ^{dim_{\C}  X} $ be the maximal unitary torus.  The blowup of $X $ along a torus-invariant submanifold $Y $ is again a toric manifold, $Bl_YX $.  As we explain in section 1.5, the action of $T $ on $X $ induces an action of $T $ on 
 $Bl_YX $.    Thus, we may study $T $-equivariant genus-0 Gromov--Witten invariants of $Bl_YX $ directly using fixed-point localization.  All faces of the  moment polytope of 
 $Y $ are faces of the  moment polytope of 
 $X $. The moment polytope of $Bl_YX $ admits a canonical inclusion into the moment polytope of $X $, for which all faces of the moment polytope of 
 $Bl_YX\setminus  \P(N_YX)  $ are contained in faces of the corresponding same dimension of the moment polytope of $X $.  Let $v_1, \dots, v_m $ be the primitive integer normal vectors to the codimension one faces of the moment polytope of a toric manifold.  Let $m_1, \dots, m_K $ be a basis of the $\Q $ vector space \[\{(a_1, \dots  a_m)   |\sum_{i=1}  ^m a_iv_i=0
 \}    , \]consisting
 of primitive integer vectors.
The toric manifold is then recovered from symplectic
reduction referred to the matrix $\m $, whose row vectors are $m_1, \dots, m_K $.
By a mirror theorem of Givental \cite{Givental_toric}   and its extensions \cite{Iritani}, a particular family of points on the Lagrangian cone of the genus-0 Gromov--Witten theory of a toric manifold is given  by an explicit formula in terms of $m_1, \dots, m_K $.  This project has its roots in the following instructive example.  Let $E $ be the total space of the projective bundle

\[\xymatrix {{\mathbb P}  (\oplus _{j=1}  ^3{\cal O}  (-a_j)  )   \ar[d]& \\
\P^2&  A:=\{[0, z_2, z_3]\subset  \P^2\}    \ar@{_{(}->}  [l]}  \]
 described by symplectic reduction with respect to the matrix  
\[ \m_X=\left( \begin{array}  {cccccc}  1& 1& 1&-a_1&-a_2&-a_3\\
0&0&0&1&1&1 \end{array}  \right).   \]
Let $[0, 0, 1]$ be the section of $E $ that maps each point 
 $x\in \P^2 $ to the point $[0, 0, 1]$ in the fiber over $x $. When $X $ is the toric bundle $E $ and $Y $ is $[0, 0, 1](A)  $ then a calculation gives \[\m_{bl_YX}  =
\left( \begin{array}  {cccccccc}  1& 1& 1&-a_1&-a_2&-a_3&0\\
0&0&0&1&1&1&0\\
 1& 0& 0&1&1&0&-1\\ \end{array}  \right).\]
In particular, \[c_1(Tbl_YX)   =``\text{the pullback of}  \ c_1(TX)   "+2P_3.\]
 In fact our main theorem arises as a generalization of this example. Here we are using the toric mirror theorems \cite{Givental_toric}, \cite{Iritani}, \cite{Brown}   as a guide to the structure of genus-0  Gromov--Witten invariants more generally (following the initial proposals of A. Elezi and A. Givental).  Elezi's work focused on projective bundles \cite{Elezi}.  In \cite{Givental_MSRI}, Givental proposed a  toric bundles generalization of Elezi's approach  using toric mirror integral representations \cite{Givental_toric}, \cite{Iritani}.  This is an ingredient in \cite{Brown}   and in the present work.

{\bf 1.7. Organization of the text.}   We recall in section 2.2  the Atiyah--Bott fixed-point localization Theorem which implies, in particular, that any
element of $H_T^*(E_{\a(A) }  )  $ is uniquely determined by  its restrictions to the $T $-fixed strata $\eps $ of $E_{\a(A) } $.  Points $\J(z)  $ on the overruled Lagrangian cone of the genus-0 $T $-equivariant   Gromov--Witten theory of 
 $E_{\a(A) } $ are certain $H $-valued formal functions, which we study in terms of their
restrictions $\{\eps^*\J(z)   \}   $.   As we recall \cite{Brown}   in section 5, the $\H_- $ projection of each of the restrictions $\eps^*\J(z)  $ consists of two types of terms.  Namely, there are terms ii) that  form simple poles expanded as $z^{-1} $ series about non-zero $H^*_T(B)  $-values of $z $.  The remaining terms i) are polynomial in $z^{-1} $ at any given order in formal variables $t, {\tilde t}  , q, {\tilde q}  , \tau , Q $.  {\bf\em  The organising principle of the text}, formulated as Theorem 2, characterizes the Lagrangian cone of the  genus-0 $T $-equivariant  Gromov--Witten theory of $E_{\a(A) } $ in terms of two conditions i)    ii)    on $\{\eps^*\J(z)   \}   $. The condition ii) says that the residues of  ${\eps}  ^*\J $ at its simple poles at non-zero values of $z $ are governed recursively with respect to $\{{\eps}  ^*\J(z)   \}   $.
The  condition i) describes the remaining poles at $z=0 $ in terms of a certain twisted Lagrangian cone of the stratum $\eps $.  The Main Theorem gives a family of points $I_{E_{\a(A) }} $  whose restrictions we check satisfy the conditions of Theorem 2.  In section 6 we verify condition ii)    for the
restrictions $\{\eps^*I_{E_{\a(A) }}  \}   $ directly, using their defining formulae. In section 7 we verify condition i)    using transformation laws \cite{Coates-Givental}   of Lagrangian cones with respect to the twisting construction from sections 1.4 and 1.5.  A new aspect of the present work relative to toric bundles is that ii)    relates the series $\{\eps^*\J(z)   \}   $ that, according to condition i), lie in Lagrangian cones derived from  genus-0 Gromov--Witten invariants of $B $ {\em and}   of $A_a $, respectively.  The Quantum Lefschetz Theorem relates the Lagrangian cone associated to the genus-0 Gromov--Witten theory of $A $ with that of $B $.  If the push-forward ${i_{A_a}}  _*:H_2(A_a, \Z)   \to H_2(B, \Z)  $ does not identify the Mori cone of $A_a $ with  that  of $B $, the opposite relation describing the  Lagrangian cone of $B $ in terms of that  of $A_a $ is realised algebraically by the Birkhoff factorization procedure and dividing by powers of $z $.  Division by $z $ does not preserve the Lagrangian cone, so we must then {\em clear denominators on both sides.}    For each $D\in MC(B)  $, denote the greatest power of $z $ that  we divide by up to order $Q^D $ in this process by $ht_A(D)  $.  We work out an example where $A $ is a smooth quintic 3-fold.

It suffices without loss of generality to assume that $\{A_a\}   $ is a single connected manifold $A $, as regards most aspects of the project.  In case there is a subtlety, we address it as it arises.

{\bf  1.8. $T $-fixed strata of $E_{\a(A) } $.}  
The $T $-fixed strata of the blowup $E_{\a(A) } $ are  in comparison  with those of $E $ as follows.  
The stratum $\a(B)  $ of $E $ is replaced by $Bl_{\a(A) }  \a(B)   =:{\tilde \a} $, which is canonically diffeomorphic to $\a(B)  $.   For each $T $-fixed section $\c\neq\a $ of $E $, the strata $\c(B)  $ of $E $ is canonically a stratum of $E_{\a(A) } $ that  we also  denote by $\c(B)  $.  Lastly, there are $T $-fixed strata of $E_{\a(A) } $ that  have no counterpart in $E $. Namely, each $T $-equivariant line bundle summand of $N_{\a(A) }  E $ gives rise to a $T $-fixed section over $A $ in the exceptional divisor.  In particular, the summand 
 $i_A^*L $ gives rise to a section 
\[ [1, 0, \dots, 0]:=\P(N_{\a(A) }  \a(B)  )   \subset  {\tilde \a}  \]
over $A $ in the total space of the exceptional divisor.  As we will see in section 2.2, the relation between $\tilde \pi^*\a^*pr_E $ and $\tilde \a^* $ is  \[\tilde \a^*-\tilde \pi^*\a^*pr_E=[1, \vec 0]^* , \]where $pr_E $ is the projection of $H^*(E_{\a(A )}  ) $ onto $H^*(E) $, and $\tilde\pi $ is the blowdown.
From now on let the symbol $\c $ stand for the $T $-fixed strata denoted $\c $ above, or for ${\tilde \a}  (B\setminus A) $.
  Let us denote the situation of a torus fixed point $\b $  connected to $\a $ by a 1-dimensional edge of the momentum polyhedron of a fiber of $E $, by $\b\sim\a $.  In this case $|\a\cup\b|=K+1 $, $|\a\cap\b|=K-1 $, and $\C^{\a\cup\b}  //T^K=:\C P^1_{\a, \b} $.  Let $j_-(\a, \b)  $ be the coordinate from $\a\setminus \a\cap\b $  and $j_+(\a, \b)  $ the coordinate from $\b\setminus \a\cap\b $.  Similarly, we have the notation $\a(\b, j_-)  $ and $\b(\a, j_+)  $.  In the next section we enhance this description of the $T $-fixed points of $E $ to a description of the $T $-fixed points of $E_{\a(A) } $.

\section{Geometry of $E_{\a(A) } $}  
{\bf 2.1. Geometric preliminaries and  decomposition of cohomology.}  
The action of $T $ on $E $ decomposes $i_{\a(A) }  ^*TE/T\a(A)  $ into a direct sum of 1-dimensional eigenspaces, 
 
\begin{align*}  N_{\a(A) }  E=
i_{\a(A) }  ^*TE/T\a(A)   \simeq i_{\a(A) }  ^*T\a(B)/T\a(A)   \oplus   i_{\a(A) }  ^*N_{\a}  E\simeq\\ i_A^*L\oplus   i_{\a(A) }  ^*N_{\a}  E.\end{align*}  
   
Let $jj_+ $ be an ordering index of these eigenspaces, where  the index value $jj_+=[1, \vec 0]$ corresponds to the bundle $i_A^*L $, and $jj_+=j_+ $ indexes the summand of $N_{\a}  E $ with $T $-weight $\a^*U_{j_+} $.  Denote the $T $-fixed section of $\P(N_{\a(A) }  E)  $ corresponding to the index $jj_+ $ by $ (\a, jj_+)  $. In the $l $-blowup case we need to include the index $a $, to denote the divisor of $B $ we blow up along.  The strata $ (\a, j_+)  $ is connected to the strata $\b(\a, j_+)  $ by the $T $-invariant edges $\C P^1_{(\a, j_+)   , \b} $.  Denote by $\chi_{\eps, \eps'}  \in H^2_T(\eps) $ the $T $-weight of $T_{\eps}  \C P^1_{\eps, \eps'} $.

Denote by $\pi:E_{\a(A) }  \to B $ the composition of ${\tilde \pi}  :E_{\a(A) }  \to E $ with the projection to the base $B $.  It is now mandatory that  we introduce the diagram

\[ \xymatrix {E_{\a(A) }  \ar@/_1pc/[dd]_{\pi}  \ar[d]^{\tilde {\pi}}  & {\mathbb P}  (N_{\a(A) }  E)   \ar@{_{(}->}  [l]&  (\a, jj_+)   (A)   \ar@{_{(}->}  [l]\\
E\ar[d]& & \a(A)   \ar@{_{(}->}  [ll]\ar@/_1pc/[u]_{(\a, jj_+) }  \\
B\ar@/_1pc/[u]_{\a}  & &  A\ar@{_{(}->}  [ll]\ar@/_1pc/[u]_{\a}}  \]
Let $N^{(\a, jj_+) } $ be the normal bundle within $E_{\a(A) } $ to the $T $-fixed section over $A $
with index $jj_+ $ in ${\mathbb P}  (N_{\a(A) }  E)  $.  

{\tt Proposition.}  {\em The action of $T $ on $E_{\a(A) } $ decomposes
 $N^{(\a, jj_+) } $ into a direct sum of $T $-equivariant line bundles, 
whose $T $-equivariant Euler classes are the elements of the set}  

\[\begin{array}  {cl}  \{-(\a, jj_+)^*U_{jj_+}  +(\a, jj_+)^*U_{j}  \}    _{j\neq jj_+, \notin\a }  \cup &\\\{-(\a, jj_+)^*U_{jj_+}  +(\a, jj_+)^*
c_1(L)\}    \cup\{(\a, jj_+)^*U_{jj_+}  \}    & , jj_+=j_+, \\

\{-(\a, jj_+)^*c_1(L)+(\a, jj_+)^*U_j\}    _{j\neq j  j_+, \notin\a }  \cup\{(\a, jj_+)^*c_1(L)   \}    & , jj_+=[1, {\vec 0}  ].\end{array}  \]

Let us now turn attention to  the restriction
map $H^*(E_{\a(A) }  )   \rightarrow H^*(\pi^{-1}  (A)  )  $.

Denote ${\tilde P} $ the $T $-equivariant Euler class of the ${\cal O}  _T(1)  $ bundle on the exceptional divisor.  
By the Lerray-Serre theorem, 
\[ H^*(\pi^{-1}  (A)  )   \simeq  H^*(A)   [{\tilde P}  ].\]
In the following section we extend the definition of $\tilde P $ to the entire blowup.  With this interpretation of $\tilde P $, recall the isomorphism of vector spaces \cite{Griffiths-Harris}  

\[    H^*(E_{\a(A) }  ) \simeq H^*(E) \oplus \left(H^*(\pi^{-1}  (A))/H^*(A)\right ), \]where the quotient is an additive quotient and $H^*(E)\simeq Im(\tilde {\pi}  ^*) $.
 On the other hand, $H^*(E)\simeq  H^*(B)[P_1, \dots, P_K]$.  The restriction of $P_i $ to the exceptional divisor is $i_{\a(A)}  ^*P_i $.

 The restriction map $i_A^* :H^*(B)   \to H^*(A)  $ and the Poincare pairing give the orthogonal  projection $\pi^{\perp} $:
 
\[\xymatrix {0\ar@{->}  [r]&Ker(i_A^*)   \ar[r]& H^*(B)   \ar@{->}  [r]\ar@{->}  [r]^{i_A^*}  &Im(i_A^*)  \ar@{->}  [r]\ar[d]_{\subset}  &0.\\
&& & H^*(A)   \ar@/_1pc/[u]_{\pi^{\perp}}  &}  \]
  
{\bf  2.2. Fixed-point localization.}   For each $\eps\in E_{\a(A) }  ^T $, the action of $T $ on $E_{\a(A) } $ decomposes $ N_{\eps}  E_{\a(A) } $ into a direct sum of 1-dimensional eigenspaces.  Denote by ${\cal N}  ^{{\eps}} $ the subbundle obtained by direct subtracting the direct summand where $T $ acts trivially.  Let $U_{A, jj}  \in H^2_T(\P (N_{\a(A) }  E))  $ be the $T $-equivariant Poincar\'e duals of the torus-invariant divisors on the exceptional divisor \linebreak ${\mathbb P}  (N_{\a(A) }  E)  $:
 \[U_{A, jj}  =\left\{\begin{array}  {ll}  {\tilde P}  +i_{\a(A)}   ^*\a^*U_j& , jj=j\notin\a \\{\tilde P}  +i_{\a( A)}  ^*c_1(N_{A}  b)   &, jj=[1, {\vec 0}  ]\end{array}  \right.        \in H^2_T(\P (N_{\a(A) }  E))   .\]

The Atiyah--Bott Theorem says that the pairing of a class $f\in H_T(E_{\a(A) }  )  $ against the fundamental class of $E_{\a(A) } $ is given by
  \begin{align*}  \int_{E_{\a(A) }}  f= \sum_{\c}  \int_{\c(B) }  \Res_{U_{j_1}  =\dots=U_{j_K}  =0}  &\frac{f dP_1\wedge\cdots\wedge dP_K} {U_1\cdots U_n}det(dU/dP)
+\\\int_{(\a, [1, {\vec 0}  ])   (A) }  \Res_{U_{A, [1, {\vec 0}  ]}  =0}  \frac{f d{\tilde P}}  {\prod_{jj'_+}  U_{A, jj'_+}}  
+\\\sum_{jj_+\neq [1, {\vec 0}  ]}  \int_{(\a, jj_+)   (A) }  &\Res_{U_{A, jj_+}  =0}  \frac{f d{\tilde P}}  {-{\tilde P}  \prod_{jj'_+}  U_{A, jj'_+}}  .\end{align*}  Namely, we sum over each of the $T $-fixed strata $\eps\in E_{\a(A) }  ^T $ the pairing of the class \[\frac{\eps^*f}  {Euler_T({\cal N}  ^{\eps}  ) }  \in H_T(\eps)   \]against the fundamental class of $\eps $, where 
\begin{align*}   Euler_T({\cal N}  ^{\eps}  )=:e^{\eps}  =&\left\{\begin{array}  {cl}  \prod_{j\notin\c}  \c ^*U_j&\eps=\c\\\prod_{jj_+\neq [1, {\vec 0}  ]}  (\a, [1, {\vec 0}  ])   ^*U_{A, jj_+}  &\eps=(\a, [1, {\vec 0}  ])   \\(\a, j_+)   ^*(-{\tilde P}  \prod_{jj'_+\neq   j_+}  U_{A, jj'_+}  )  &\eps=(\a, j_+)   .\end{array}  \right.       \end{align*}  

Thus, denote ${\tilde P} $ the class in $H_T^*(E_{\a(A)}  , \Q(\l)  )  $ that  restricts to the $T $-equivariant Euler class of the ${\cal O}  _T(1)  $ bundle on the exceptional divisor, and 
restricts to zero at all $T $-fixed strata in the complement of the exceptional divisor.  This means in particular, that $(\a, [1, {\vec 0}  ])   ^*{\tilde P}  =-i_A^*c_1(L)  $ and $({\tilde \pi   ^*\a}  ^*pr_E){\tilde P}  =0 $.

In particular, \[c_1(TE_{\a(A) }  )   ={\tilde \pi}  ^*c_1(TE)+(N-K)   \tilde P.\]

Define a $T $-equivariant line bundle ${\tilde l} $ 
 over the union of torus-invariant edges of 
 $E_{\a(A) } $ as follows.  It restricts to the ${\cal O}  _T(1)  $ bundle over the edges of the exceptional divisor, restricts to the trivial bundle  over the edges $\C P^1_{\c, \c'} $ and whose $T $-equivariant Euler class restricts to $c_1^T(T^*\C P^1_{(\a, j_+)   , \b}  )+{\tilde \pi}  ^*U_{j_-(\a, \b)}   $ over the edges $\C P^1_{(\a, j_+)   , \b} $.

{\tt Proposition.}  {\em The ${\tilde P} $ pairings on elements of $H_2(E_{\a(A) }  , \Z)  $ take values in $\Z $.}  

{\em Proof.}   The restriction of ${\tilde P} $ to the union of torus invariant edges coincides with the class $c_1^T({\tilde l}  ) $.
Apply the Atiyah--Bott fixed-point localization Theorem to the restriction of ${\tilde P} $ to the union of torus-invariant edges of
 $E_{\a(A) } $, 
 \[{\tilde P}  (d_{(\a, j_+)   , \b}  )   =\frac{ (\a, j_+)   ^*{\tilde P}  -0}  {\chi_{(\a, j_+)   , \b}  
}  =-1, \]and ${\tilde P}  (d_{\a, \b}  )   =0 $ for all $\b\sim \a $. 
Thus, ${\tilde P} $ induces an element of $H^2(E_{\a(A) }  , \Z)  $.

{\tt Proposition.}  {\em
Denote $i:\C P^1_{\a, \b}  \to E_{\a(A) } $ the 
inclusion map.  The normal subbundle
 \[\NN^{\tilde \a(B\setminus A)}  ={\tilde \a}  (B\setminus A)^*TE_{\a(A) }  /T{\tilde \a}  (B\setminus A)  \]  extends over $\C P^1_{\a, \b} $ as $T\C P^1_{\a, \b}  \oplus _{j\notin\a \cup\b}  l_j $, where $l_j $ are $T $-equivariant line bundles with
first Chern classes $i^*{\tilde \pi}  ^*U_j $. 
  
Denote $i:\C P^1_{(\a, j_+)   , \b}  \to E_{\a(A) } $ the 
inclusion map.  The normal bundle
\[\NN^{(\a, j_+) }  =(\a, j_+)^* TE_{\a(A) }  /T(\a, j_+)\]extends over $\C P^1_{(\a, j_+)   , \b} $ as $T\C P^1_{(\a, j_+)   , 
\b}  \oplus _{j\notin\a \cup\b}  l_j $, 
 where $l_j $ are $T $-equivariant line bundles with
first Chern classes $i^*{\tilde \pi}  ^*U_j+(c_1^T(T^*\C P^1_{(\a, j_+)   , \b }  )+i^*{\tilde \pi}  ^*U_{j_-(\a, \b)}  ) $.

Denote $i:\C P^1_{(\a, j_+)   , (\a, jj_+') }  \to E_{\a(A) } $ the 
inclusion map.  The normal bundle
\[\NN^{(\a, j_+) }  =(\a, j_+)^* TE_{\a(A) }  /T(\a, j_+)\]extends over $\C P^1_{(\a, j_+)   , (\a, jj_+') } $ as $T\C P^1_{(\a, j_+)   , 
(\a, jj_+')   }  \oplus _{j\notin \a, \neq j_+, jj_+'}  l_j $, 
 where $l_j $ are $T $-equivariant line bundles with
first Chern classes $i^*{\tilde \pi}  ^*U_j+i^*{\tilde P} $, while the normal bundle
\[\NN^{(\a, [1, \vec 0]) }  = (\a, [1, \vec 0])^*(TE_{\a(A) }  /T(\a, [1, \vec 0]))   \ominus i^*_AL\]extends over $\C P^1_{(\a, [1, \vec 0])   , (\a, j_+') } $ as $T\C P^1_{(\a, [1, \vec 0])   , (\a, j_+') }  
 \oplus _{j\notin\a , \neq j_+'}  l_j, $
where $l_j $ are $T $-equivariant line bundles with
first Chern classes $i^*{\tilde \pi}  ^*U_j+i^*{\tilde P} $.}

\section{The $ht_A $ function}

Let $P_{1}  , \dots, P_{r} $ be a basis of $H^{1, 1}  (B, \Z)  $, and $P_{r+1}  , \dots, P_{r+2s} $  a basis of $H^{0, 2}  (B, \Z)\linebreak  \oplus   H^{2, 0}  (B, \Z)  $, with dual bases $\tau_{1}  , \dots, \tau_{r} $ and $\tau_{r+1}  , \dots, \tau_{r+2s} $.  Let $\vartheta_1, \dots, \vartheta_v $ be coordinates on $H^*(B)/<H^2(B)   , +> $.  Define \[\J_{\{e(\cdot)   , L_a\}}  (-z, \vartheta+\tau)   := \sum_{D\in MC(B) }  Q^D\J^D_B(-z, \vartheta+ \tau)\prod_{a=1}  ^l\prod_{m=1}  ^{c_1(L_a)   (D) }  (c_1(L_a)   -mz). \]
{\tt Quantum Lefschetz Theorem \cite{Coates-Givental}.}    {\em Suppose $c_1(L)   (D)   \geq 0\ \forall D\in MC(B)  $, or more generally that $L $ is convex.  
Then for each $\vartheta+ \tau\in H^*(B)  $ and for each smooth family $\J_B(-z, \tau)   \subset  \LL_B $, the series  modification 
\begin{align*}  i_A^*{\tilde \J}  _{e(\cdot)   , L}  (-z, \vartheta+\tau)   :=& \sum_{D\in Im(i_A)   _*\subset  MC(B)   }  Q^D\times \\
&i_A^*\left(\J^D_B(-z, \vartheta+ \tau)\prod_{m=1}  ^{c_1(L)   (D) }  (c_1(L)   -mz)    \right)   \end{align*}  lies in the image by $\pi^{\perp} $ of the Lagrangian cone associated to the genus-0 Gromov--Witten theory of $A $ with inputs from $Im(i_A^*)  $.}  

Let us assume that $\J_B(z, \vartheta+\tau)  $ has the property (Div+Str primary)   that its dependence on $\tau_{1}  , \dots, \tau_{r+2s} $ is of the form \[\sum_{D\in MC(B) }  a_D(z, \vartheta)  e^{\tau/z}  Q^De^{\tau(D) }  , \]
where $a_D(z, \vartheta)  $ do not depend on $\tau_{1}  , \dots, \tau_{r+2s} $, are Laurent polynomials \linebreak in $z $ valued in $H^*(B, \Q(\l)  )  $, and $a_0(z, \vartheta)   =e^{\vartheta/z} $. 
Then, both series \linebreak $i_A^*\J_{e(\ \cdot\ )   , L}  (-z, \vartheta+ \tau)=: ze^{i_ A^*P_Bt/z}  (1+B_1q+B_2q^2+\dots) $ and $i_A^*{\tilde \J}  _{e(\ \cdot), L}  (-z, \vartheta+\tau)=: ze^{i_ A^*P_Bt/z}  (1+A_1q^5+A_2q^{10}  +\dots) $ have the property Div+Str primary.
Also assume that \[i_ A^*e^{-(\vartheta+\tau)   /z}  z^{-1}i_A^*{\tilde \J}  _{e(\ \cdot), L}  (-z, \vartheta+\tau)   =1+{\cal O}  (Q), \]
which is equivalent to the condition 
	 \[i_ A^*e^{-(\vartheta+\tau)   /z}  z^{-1}   \J_{e(\cdot), L}  (z, \vartheta+\tau)   =1+{\cal O}  (Q).\]
Our goal is to prove well-definedness of the least positive integer function $ht_A:MC(B)   \to \N $ such that, for each $D\in MC(B)  $, the truncation of  
 $z^{ht_A(D) }   i_A^*\J_{e(\ \cdot\ )   , L}  (z, \vartheta+ \tau)   \mod (Q^D) $ is a formal linear combination of vectors in the linear space 
\begin{align*}   zT_{\f(-z) }  (\pi^{\perp}  \LL_ A)   \subset  \pi^{\perp}  \LL_ A\mod (Q^D)   , \ \ \text{where}       \\ \f(z)   := z   ^{ht_A(D) }i_A^*{\tilde \J}  _{e(\ \cdot), L}  (-z, \vartheta+\tau)  .\end{align*}  We now prove well-definedness of $ht_A $ by giving a combinatorial algorithm for computing it.
  We observe the following (Divisor-, String-)    differential equations
\begin{align*}  i_ A^*P_ie^{(\vartheta+\tau)/z}  &Q^De^{\tau((i_A)_*D) }  = \\&\left\{\begin{array}  {cl}  (z\p_{\tau_i}  -z
P_i(i_A)_*D  ) i_ A^* e^{(\vartheta+\tau)/z}  Q^De^{\tau(i_A)_*D }  &i=1, \dots, r\\
 z\p_{\tau_i}     i_ A^*e^{(\vartheta+\tau)/z}  Q^De^{\tau(i_A)_*D }  &i=r+1, \dots, r+2s\end{array}  \right.       \\&=:({\vec P}  _i)_D i_A^*e^{(\vartheta+\tau)/z}  Q^De^{\tau(i_A)_*D }  i=1, \dots, r+2s.\end{align*}  
For any polynomial $\phi $ in variables $P_1, \dots, P_{r+2s}  , z $ with coefficients in $\QQ $, it follows that 
\begin{align*}  i_ A^*\phi(P_1, \dots, P_{r+2s}  , z)   &e^{(\vartheta+\tau)/z}  Q^De^{\tau(i_A)_*D }  =\\&\phi(({\vec P}  _1)_D, \dots, ({\vec P}  _{r+2s}  )_D, z) i_ A^* e^{(\vartheta+\tau)/z}  Q^De^{\tau(i_A)_*D }  .\end{align*}  

	Identify the Mori cone of the base with $\N^r $, and decompose $\N^r $ into hyperplanes:  
\[ \N^r =\coprod_{c=1}  ^{\infty}  \Delta_c, \]
where\[ \Delta_c:=\{(a_1, \dots, a_r)   \in\N^r| \sum_{i=1}  ^ra_i=c\}.\]
Let $\iota:\N^r\to \N $ be any  bijective map such that  for each $c<c' $, $x\in\Delta_c $ and $y\in \Delta_{c'} $, it follows $\iota(x)   <\iota(y)  $. 

Define $\{C_n\}   $ recursively:\[ ( 1+B_1q+B_2q^2+\dots)=\\(1+C_1q+C_2q^2+\dots)   (1+A_1q^5+A_2q^{10}  +\dots)   .\]
Define $\phi_n $ through the formula

\begin{align*}  i_A^*\J_{e(\cdot)   , L}  ^{\iota^{-1}  (n) }  (z, \vartheta+\tau)  &  =\\\sum_{n=0}  ^ {\infty}  Q^{\iota^{-1}  (n)}  &\frac{\phi_n(({\vec P}  _1)_{\iota^{-1}  (n)}  , \dots, ({\vec P}  _{r+2s}  )_{\iota^{-1}  (n)}  , z)}  {z^{d_n}}  
e^{\tau(\iota^{-1}  (n)  )} i_A^*{\tilde \J}  _{e(\ \cdot), L}  (-z, \vartheta+\tau), \end{align*}      
where $d_n $ is the (maximal) pole order of $C_n $ at $z=0 $, so that 

\begin{align*}   &C_n=\frac{e^{\tau(\iota^{-1}  (n)  ) }}  {z^{d_n}}  \times \\&\left(\phi_n(i_A^*P_1, \dots, i_A^*P_{r+2s}  , z)+\phi_{n-1}  (i_A^*P_1-zP_1(i_A)_*\iota^{-1}  (n-1), \dots, z)+\cdots+\right.       \\&\left.\phi_1(i_A^*P_1-zP_1(i_A)_*\iota^{-1}  (1), \dots, z)\right).\end{align*}  

Define \[ht_A(\iota^{-1}  (n)  )   :=max\ \{\{d_{n'}  \}    _{n'\leq n}  , 0\}    , \]and \[ht_{\{A_a\}  }  (\iota^{-1}  (n)  )   :=max_a\ ht_{A_a}  (\iota^{-1}  (n)  )   .\]
Let $\J_{A_a}(-z, \vartheta+\tau) $ be the unique family of points of $\LL_{A_a} $ for which $\pi^{\perp}\J_{A_a}(-z, \vartheta+\tau)=(-z)^{ht_{\{A_a \}}  (D') } i_{A_a}^* \J_{e(\ \cdot), L_a}(-z, \vartheta+\tau) $.

\[\G_{A_a}  (-z, \vartheta+\tau):=\J_{A_a}(-z, \vartheta+\tau)-\pi^{\perp}\J_{A_a}(-z, \vartheta+\tau). \]

\[\G^D_{\{A_a\}}  (-z, \vartheta+\tau):=\oplus_{a=1}^l\G^D_{A_a}(-z, \vartheta+\tau)\prod_{a'=1, \neq a}^l\prod_{m=1}^{c_1(L_{a'})   (D)   }(c_1(L_{a'})   -mz) . \]

{\bf Example.}   Let $b=\C P^4, L={\cal O}  (5)   \to \C P^4 $, and $A\subset   \C P^4 $  a smooth quintic 3-fold.

 Let $P_B $ be the K\"ahler generator of $H^2(B, \Z)  $, and take $\J_B(z, \tau)  $ to be the J-function of $\C P^4 $ at the point $\tau=P_Bt $, 
\[\J_B(z, \tau)   =ze^{P_Bt/z}  \sum_{d=0}  (qe^t)   ^d\frac{1}  {\prod_{m=1}  ^d(P_B+mz)   ^5}  , \]
\begin{align*}  i_A^*\J_{e(\ \cdot)   , L}  (z, \tau)   =i_A^*ze^{P_Bt/z}  \sum_{d=0}  (qe^t)   ^d\frac{\prod_{m=1}  ^{5d}  (5P_B+mz) }  {\prod_{m=1}  ^d(P_B+mz)   ^5}  \text{mod}  \ P_B^4=:\\ ze^{i_ A^*P_Bt/z}  (1+B_1q+B_2q^2+\dots)\ , \end{align*}

\[{\tilde \J}  _{e(\ \cdot)   , L}  (z, \tau)   =ze^{P_Bt/z}  \sum_{d=0}  (qe^t)   ^{5d}  \frac{\prod_{m=1}  ^{25d}  (5P_B+mz) }  {\prod_{m=1}  ^{5d}  (P_B+mz)   ^5}  , \]

\begin{align*}i_A^*{\tilde \J}  _{e(\ \cdot)   , L}  (z, \tau)   =i_ A^*ze^{P_Bt/z}  \sum_{d=0}  (qe^t)   ^{5d}  \frac{\prod_{m=1}  ^{25d}  (5P_B+mz) }  {\prod_{m=1}  ^{5d}  (P_B+mz)   ^5}  \text{mod}  \ P_B^4=:\\ ze^{i_ A^*P_Bt/z}  (1+A_1q^5+A_2q^{10}  +\dots)\ .\end{align*}    Thus, we deduce the relation $A_n=B_{5n} $ mod $P_B^4, \ \forall n\geq 1 $.
The coefficient of
 $ze^{P_Bt/z}  (qe^t)  $ in $\J_{e(\ \cdot)   , L}  (z, \tau)  $ is \[\prod_{m=1}  ^5\frac{(5P_B+mz) }  {(P_B+z) }  =5\prod_{j=1}  ^4(5-j[1-(P_Bz^{-1}  )+(P_Bz^{-1}  )   ^2-(P_Bz^{-1}  )   ^3+(P_Bz^{-1}  )^4])   .\]
This latter series is easily expanded by the formal identity
\begin{align*}  5\prod_{j=1}  ^4(5-x_j)   
=5(5^4-5^3&\sum_{j=1}  ^4x_j+5^2\sum_{\text{distinct}  \ j_1, j_2}  x_{j_1}  x_{j_2}  
-\\&5\sum_{\text{distinct}  \ j_1, j_2, j_3}   x_{j_1}  x_{j_2}  x_{j_1}  +x_1x_2 x_3x_4)   .\end{align*}  Thus, for example, the coefficient of $(P_Bz^{-1}  )   ^0 $ of the series 
 $\prod_{m=1}  ^4\frac{(5P_B+mz) }  {(P_B+z) } $ is 
\[5\sum_{k=0}  ^4(-1)   ^k 5^{4-k}  {\o {Sym}}  _k(1, 2, 3, 4)   .\]where ${\o {Sym}}  _k(1, 2, 3, 4)   :=\sum_{\text{distinct}  \ j_1, \dots , j_K \in\{1, 2, 3, 4\}  }    j_1\cdots j_K $.
The remaining powers of $P_Bz^{-1} $ are easily restored, 
\begin{align*} e^{-t}  B_1= \prod_{m=1}  ^5\frac{(5P_B+mz) }  {(P_B+z) }  \ \text{mod}  \  P_B^4=\\ 5\sum_{m=0}  ^4(-P_Bz^{-1}  )   ^m \sum_{k=0}  ^4(-1)   ^k 5^{4-k}  {\o {Sym}}  _k(1, 2, 3, 4)\times  \\ \sum_{
 (c_1\geq\dots \geq c_k)   \vdash m}  |Orbit_{S_k}  (c_1, \dots, c_k)   |\ \text{mod}  \  P_B^4.
\end{align*}  \begin{align*} e^{-t}  \frac{B_n}  {B_{n-1}}  =\prod_{m=5(n-1)+1}  ^{5n}  \frac{(5P_B+mz) }  {(P_B+nz) }  \ \text{mod}  \  P_B^4=\\ 5\prod_{j=1}  ^4(5-\frac{j}  {n}  [1-(\frac{P_Bz^{-1}}  {n}  )+(\frac{P_Bz^{-1}}  {n}  )   ^2-(\frac{P_Bz^{-1}}  {n}  )   ^3+(\frac{P_Bz^{-1}}  {n}  )^4])\ \text{mod}  \  P_B^4=\\
5\sum_{m=0}  ^3(-\frac{P_Bz^{-1}}  {n}  )   ^m \sum_{k=0}  ^4(-1)   ^k
 5^{4-k}  {\o {Sym}}  _k(\frac{1}  {n}  , \frac{2}  {n}  , \frac{3}  {n}  , \frac{4}  {n}  )\times \\   \sum_{
(c_1\geq\dots \geq c_k)   \vdash m}  |Orbit_{S_k}  (c_1, \dots, c_k)   |.\end{align*}   The latter formulae  determine $\{A_n\}   $ and $\{B_n\}   $ respectively.  Then, 
\[(1+C_1q+C_2q^2+\dots)   (1+A_1q^5+A_2q^{10}  +\dots)   =( 1+B_1q+B_2q^2+\dots)\] determines $\{C_n\}   $ recursively:
 \[C_n=B_n-C_{n-5}  A_1-C_{n-10}  A_2-\dots.\]For each $n\geq 1 $, $d_n $ is the maximal power of $P_B $ in the $B_n $ series.  Thus, $d_n=3 $ for all $n \geq 1 $.

{\tt Proposition.}   {\em If $B $ is $\C P^n $, if $A $ is the zero locus of a generic section of a convex line bundle $L $ over $B $, if $\J_B $ is the J-function of $B $, and if the class $(c_1(TB)-c_1(L)) $ of the base is nonnegative as a functional on $MC(B) $, then $ht_A(D)=(c_1(TB)-c_1(L))\cap D+(dim_{\C}  B-1)\ \forall D\in MC(B)\setminus 0 $.}  

\section{Main results}  
{\bf 4.1. The I-function.}   Upon extension of scalars $\Z\subset \Q $ of homology groups, the Mori cone of $E_{\a(A) } $  includes into $H_2(E_{\a(A) }  , \Q)  $. 
 Given ${\tilde \D}  \in H_2(E_{\a(A) }  , \Q)  $ or ${\tilde \D}  \in MC(E_{\a(A) }  )  $, define $\D:={\tilde \pi}  _*({\tilde \D}  )   , D:=\pi_*({\tilde \D}  )   , d_i:=P_i({\tilde \D}  )   , {\tilde d}_a:={\tilde P}_a({\tilde \D}  )  $, and these values uniquely determine ${\tilde \D} $.

Henceforth we use the Gamma-function convention:
                        \[ \prod_{m=1}  ^n (U+mz)    := \prod_{m=-\infty}  ^n (U+mz)   \ \ / \prod_{m=-\infty}  ^0
                        (U+mz)   .\]

{\tt Main Theorem.}   {\em Let $E $ be a toric fibration over base $B $, whose fibers are not copies of the point, and let $\a:B\to E $ be a $T $-fixed section. 
Let $L_a $ be convex line bundles over $B $, and $A_a $ smooth divisors of $B $ arising as the zero loci of generic sections of $L_a $.  Further assume the $\{A_a\}   $ to be  mutually disjoint.
{\tt  Case 1:}   If the push-forwards ${i_{A_a}}  _*:H_2(A_a, \Z)   \to H_2(B, \Z)  $ do not identify the Mori cone of $A_a $ with that  of $B $, 
 then for each $D'\in MC(B)  $, for each $(t, {\tilde t}  )  $, for each $\tau\in H^*(B)$ and for each smooth family $\J_B(-z, \tau)   \subset  \LL_B $ with the property Div+Str primary, the truncation mod $(Q^{D'}  )  $ of the $z\to -z $ version of the series $I_{E_{\coprod_a A_a}}  (z, t, \tilde t, \tau, q, \tilde q, Q)   \subset (\text{a completion}   \footnote{by an extension of the Novikov ring of $E_{\a(A)} $. See the first example of the Main theorem and the second Remark in section 6.}  \ \text{of}  )\ \H $ defined by}

\begin{align*}  I_{E_{\coprod_a A_a}}  (z, t, \tilde t, \tau, q, \tilde q, Q)   =e^{Pt/z} e^{{\tilde P}  {\tilde t}  /z}   \sum_{ d\in \Z^K, 
\tilde d\in \Z^l, D\in MC(B) }  \\\frac{(qe^t)   ^d({\tilde q} e^{\tilde t}  )   ^{\tilde d}   Q^D({\bf z^{ht_{\{A_a\}  }  (D') }}  \J^D_{\{e(\cdot), L_a\}}  (z, \tau) +\G^D_{\{ A_a\}}  (z, \tau))}  {\prod_{j\notin \{j_+(\a, \b)   |  \a\sim\b\}  }  
\prod_{m=1}  ^{U_j(\D) }  (U_j+mz) }  \times \\ \frac{1}  
{ \prod_{j\in \{j_+(\a, \b)   |  \a\sim\b\}  }  
\prod_{m=1}  ^{U_j(\D)+\sum_{a=1}  ^l{\tilde d}_a }  (U_j+\sum_{a=1}  ^l{\tilde P}_a+mz) }  \times \\
\frac{1}  {{\bf \prod_{a=1}  ^{\ell}  \prod_{m=1}  ^{c_1(L_a)   (D)+{\tilde d}_a}  (c_1(L_a)+{\tilde P}_a+mz) }}  \times \\ 
\frac{1}  {\prod_{a=1}  ^l
\prod_{m=1}  ^{-{\tilde d}_a}  (-{\tilde P}_a+mz) }  \end{align*}  {\em lies in the truncation mod $(Q^{D'}  )  $ of the Lagrangian cone associated to the genus-0 Gromov--Witten theory of the blowup within $E $ of $\a(\coprod_a A_a)  $.
{\tt  Case 2:}   If the push-forwards ${i_{A_a}}  _*:H_2(A_a, \Z)   \to H_2(B, \Z)  $ identify the Mori cone of $A_a $ with that  of $B  $, then $ht_{\{A_a\}  }  (D')   =0\forall D'\in MC(B)  $, and the preceding series  lies in the preceding cone without any truncation condition on either, and without assuming the property Div+Str primary for the smooth family $\J_B(-z, \tau)   \subset  \LL_B $.}  

{\tt Remark.}    When the fibers are copies of the point then we omit the sum over ${\tilde d} $ and we set ${\tilde P} $ to zero, since the projective fibers are also  copies of the point.
Keeping these interpretations in mind, the theorem remains true when the fiber of the toric fibration is the point.  The theorem reduces to the 
statement $\J_B\subset  \LL_B $.

{\tt Remark.}    The natural generalization of the Main Theorem to the case of several $T $-fixed sections of $E $ coincides, at a first level of analysis, with the natural generalization of the mirror theory of section 7.

{\tt Remark.} The proof of Theorem 2 indicates the dependence of points of $\LL_E$ upon domain variables  from $ H^*(A)$. 

{\tt Conjecture.} The dependence on domain variables \[ (u,\dots,u)\in \oplus_{jj_+} H^*((\a,jj_+))/H^*(A)\]may be incorporated into the Main Theorem by replacing $\tau\to\tau +\eps _{\J}(\tau,u)$ in the argument of $\J_B(\tau,z)$ and $ i^*_A\tau\to i^*_A\tau+u$ in the argument of  $\G_A(\tau,z)$, for some function $\eps _{\J}:H^*(B)\oplus H^*(A)\rightarrow H^*(B)$, $\eps _{\J}(\tau,0)=0$.

{\tt Some examples of the main Theorem.}  

1. Let $B $ be a smooth toric variety obtained by $T^K $-symplectic reduction of $\C^N $ and $A $ a (nef)    coordinate hyperplane divisor of $B $.
An instance of $E_{\a(A)} $ in this case is the example in section 1.6.  If the bundle $L $ is considered as $T^3 $-equivariant, then $\tilde P $ is $T^6 $-equivariant.
The class $P_3\in H_{T^7}  ^2(E_{\a(A)}  ) $ is not the same equivariantly as $\tilde P\in H_{T^6}  ^2 (E_{\a(A)}  ) $.  Modulo this difference, the series of the Main Theorem is an extension outside the Novikov ring of the series of the toric mirror theorems, in the same example and generally for symplectic toric manifolds \cite{Givental_toric}, \cite{Iritani}, \cite{Brown}.  

2. Let $B $ be $\P^2\times (\P^2)   ^* $, $L={\cal O}  _{\P^2}  (1) \otimes {\cal O}  _{(\P^2)^*}  (1) $  and $A $ the manifold of complete flags in $\C^3 $.

{\tt  Corollary.}  {\em Let $E $ be a toric fibration over base $B $, whose fibers are not copies of the point, 
 and let $\a:B\to E $ be a $T $-fixed section.  Then for each $(t, {\tilde t}  )  $, for each $\tau\in H^2(B)   $ and for each smooth family $\J_B(-z, \tau)   \subset  \LL_B $, the $z\to -z $ version of the series}  
\begin{align*}  I_{E_{\a}}  (z, t, \tilde t, \tau, q, \tilde q, Q)   =e^{Pt/z} e^{{\tilde  P}  {\tilde  t}  /z}   \sum_{ d\in \Z^K, 
\tilde d\in \Z, D\in MC(B) }  \\\frac{(qe^t)   ^d({\tilde q} e^{\tilde t}  )   ^{\tilde d}  Q^D\J^D_B(z, \tau) }  {\prod_{j\notin \{j_+(\a, \b)   |  \a\sim\b\}  }  
\prod_{m=1}  ^{U_j(\D) }  (U_j+mz) }  \times \\ \frac{1}  
{ \prod_{j\in \{j_+(\a, \b)   |  \a\sim\b\}  }  
\prod_{m=1}  ^{U_j(\D)+{\tilde d}}  (U_j+{\tilde P}  +mz)   
}  \times \\ 
\frac{1}  {\prod_{m=1}  ^{-{\tilde d}}  (-{\tilde P}  +mz) }  \end{align*}  {\em lies in the Lagrangian cone associated to the genus-0 Gromov--Witten theory of  the  blowup within $E $ of $\a(B)  $.}

{\tt Application to codimension $>1 $ subvarieties $A\subset  B $.}  
 Let $E\to B $ be a symplectic reduction of a direct sum of line bundles pulled back from $A $, and $B\to A $ a symplectic reduction of a direct sum of line bundles also pulled back from $A $.
 $T $-fixed sections  $\a_1: A\to B $ and $\a_2:B\to E $ may be considered as index subsets, respectively.  The disjoint union of index subsets defines a $T $-fixed section $\a_1\coprod\a_2:A\to E $.  Then Corollary applies to $Bl_{\a_1\coprod\a_2(A) }  E $, where  the matrix used for the symplectic reduction is block diagonal with a block for each of the fibers. 
  
\section{The $T $-equivariant cone $\LL_{E_{\a(A) }} $}  
{\bf  5.1.  Localization of stable maps.}    The work of Graber--Pandharipande 
\cite{Graber-Pandharipande}    justifies the fixed-point localisation technique for computing integrals of $T $-equivariant cohomology classes over virtual   fundamental cycles in the moduli spaces of stable maps to $E_{\a(A) } $.  Here the $T $-equivariant normal ``bundle" to a $T $-fixed stable map is actually a virtual   (orbi-)    {\em virtual  bundle}   in $T $-equivariant K-theory.  The connected components of the $T $-fixed loci in the moduli spaces of genus-0 stable maps are fiber products of moduli spaces of genus-0  stable maps into the $T $-fixed strata of $E_{\a(A) } $.  
Let $C $ be a {\em leg}   of $ \Sigma $; i.e., an irreducible component of $\Sigma $ that maps surjectively to a $T $-invariant edge of $E_{\a(A)} $.
The fiber product is defined by reference to the curves from $\eps_{0, n, D} $, 
 from $\eps'_{0, n, D} $ and toric edges $f(C) $.  The image points $f(0) $ and $f(\infty) $ coincide with the images of the marked points of  stable maps from $\eps_{0, n, D} $ and from $\eps'_{0, n, D} $ in their roles as nodal points.  There is also the case that either $0\in C $ or $\infty\in C $ may be a marked or unmarked point of $\Sigma $, not connecting $C $ to any other curve component of $ \Sigma $.

There are three disjoint cases to consider, depending upon how the 1-dimensional $T^{\C} $ orbit $f(C)  $ intersects the exceptional divisor.  Equivalently, these cases are distinguished by the blowdown image of the point set $f(C) $.   Firstly (2.bb), the projection of the toric edge along the blowdown map is again a toric edge at each point of the given fiber product.
 Suppose that  the two strata connecting a toric edge map via the blowdown to the $T $-fixed sections $\a $ and $\b $.  The factor of the fibre product given by genus-0 stable maps into ${\tilde \a}  :=Bl_{\a(A) }  \a(B)   \simeq \a(B)  $ can be non-compact, as follows.  Given a toric edge connecting 
 $Bl_{\a(A) }  \a(B)   \setminus  \P(N_{\a(A) }  \a(B)  )  $ to $\b(B\setminus  A)  $, the nodal point in $\a(B\setminus  A)  $ is unable to enter the exceptional divisor.
 The Atiyah--Bott formula implies that the correct cohomology group to use for the non-compact space 
 $Bl_{\a(A) }  \a(B)    \setminus  \P(N_{\a(A) }  \a(B)  )  $ is the pullback to $H^*(Bl_{\a(A) }  \a(B)    \setminus  \P(N_{\a(A) }  \a(B)  )  )  $ of $H^*(B)  $.  A similar case to consider is when the toric edge connects to $\c, \c'\simeq B$, where $\c\neq \c' $. 

Secondly (2.ab), the torus-fixed points of the toric edge connect to the rest of the $T $-fixed stable map at $\b(A)  $ and at $(\a, j_+)   (A)  $.  In this case too, the  blowdown image of the toric edge is also a toric edge.

Third (2.aa), the toric edge is contracted by the blowdown map at each point of the given fiber product. 
The torus-fixed points of the toric edge connect to the rest of the $T $-fixed stable map at $(\a, jj_+)   (A)  $ and at $(\a, jj'_+)   (A)  $.

There are three types of terms that contribute to the series ${\eps}  ^*\J(z, t)  $. Namely, the polynomial term $\eps^*t(z)   -z $, and then two types of contributions to the $\H_- $ projection of the series ${\eps}  ^*\J(z, t)  $.   Given a $T $-fixed stable map to $E_{\a(A) } $, which we denote by $(f; \Sigma, p(\Sigma)  )  $, let 
 $C $ be the smooth irreducible component of $\Sigma $ that contains the first marked point of the source of the stable map.  In order for the stable map $[f;\Sigma, p(\Sigma)   ]$ to contribute to $\eps^*\J(z, t)  $, $f $ must map the first marked point into the stratum $\eps $. The latter two types of contributions are determined by
 whether 

i)    All points of $C $ are mapped by $f $ into the $T $-fixed stratum $\eps $.  In this case let $C' $ be the maximal connected subset  of $\Sigma $ containing $C $ that maps to $\eps $, and let $D'=f_*[C']\in H_2(\eps, \Z)  $.

ii)   $C $ maps to a $T $-invariant $\C P^1 $ in $E_{\a (A) } $ connecting $T $-fixed strata $\eps $ and $\eps' $.  Let us assume that, 
in the normalization of $\Sigma $, $C $ is a $\C P^1 $ with two marked
points---which we may take to be $0 $ and $\infty $ via the action of
 $PSL_2 (\C)  $ on $\C P^1 $---, that  there is a marked point of $\Sigma $
at $0\in C $, and that  the marked point at $\infty $ corresponds to a
node of $\Sigma $.  Thus the stable map takes $C $ to a $\C
P^1_{\eps, \eps'} $, maps the first marked point of $\Sigma $ at $0\in C $ to $\eps $ and maps 
 $\infty $ to a nodal point of the stable map at $\eps' $, and as it follows from the work of Kontsevich  \cite{Kontsevich}, is given by $f([z, w])   =[z^k, w^k]\in\C
P^1_{\eps, \eps'} $.

  The group of automorphisms of $C $ that  fix $0 $ and
 $\infty $ is $\C^*: \Aut(C;0, \infty)   \simeq\C^*. $

The equivalence class of $(f;C, 0, \infty)  $ is $T $-fixed; i.e., there is 
a 1-dimensional linear representation \[ \varphi: T:\rightarrow \Aut
(C;0, \infty)   \]such that  for all $h\in T $, \[
hf([z, w])   =f([\varphi(h)   z, w])   =[\varphi^k(h)   z^k, w^k]\in\C P^1_{\eps, \eps'}  .\]

 The weight of the $T $-representation $\varphi $ is $\frac{\chi_{\eps, \eps'}}  {k} $, since the $T $ weight on the LHS is $\chi_{\eps, \eps'} $.  The map $\varphi $ induces a representation of $T $ on $T_0C $, and therefore the $T $-weight of the cotangent line bundle at the 1st marked point of $\Sigma $ (i.e., of $0\in C) $ is $-\frac{\chi_{\eps, \eps'}}  {k} $. 
The virtual  normal bundle \[ H^0(\Sigma, f^*\NN^{\eps}  )   \ominus
H^1(\Sigma, f^*\NN^{\eps}  )   \ominus (\Lie \Aut(\Sigma; p(\Sigma)  )   \cap H^0(\Sigma, 
f^*\NN^{\eps}  )  )   \]to the $T $-fixed stable map with source $\Sigma $ decomposes
into :

(i)     The virtual  normal bundle over the stable map with source
 $\Sigma':=\o{\Sigma\setminus  C} $, and

(ii)    A virtual  {\em vector space}  

\[ {\mathcal {N}}  _{\eps, \eps'}  (k)   \simeq H^0(C, f^*\NN^{\eps}  )   \ominus
H^1(C, f^*\NN^{\eps}  )   \ominus \Lie \Aut(C; 0, \infty)   
\ominus \NN^{\eps'}  \]
over the {\em point}  
 $[f;C, 0, \infty]$.  This virtual  vector space is the fiber of a virtual  bundle. We use the 
same notation for the bundle as for the fiber.

There is a subtlety to address here.  Namely, we first meet $H^1(\Sigma, f^* l_j) $ in the  deformation theory description of the virtual normal bundle. The description of this cohomology group in terms of $H^1(C, f^*l_j) $ and $H^1(\Sigma', f^*l_j) $ requires some work.
   Consider a contractible
neighborhood of the node in $\Sigma $ where $C $ meets $\Sigma' $.
Recall that $w $ is a local coordinate in $C $ centered at the node, and
denote by $z' $ a local coordinate on $\Sigma' $ centered at the node.
Then, the local coordinate expression
\[ \frac{dw\wedge dz'}  {d(wz')}  =\frac{dw}  {w}  -\frac{dz'}  {z'}   \]
defines a regular non-vanishing 1-form in the contractible
neighborhood of the node.  We pull this 1-form back along the map that
normalizes the node, and denote $\omega_C $ the restriction to $C $ and
 $\omega_{\Sigma'} $ the restriction to $\Sigma' $.  Observe that
 $\omega_C $ has a simple pole at $w=0 $, and $\omega_{\Sigma'} $ has a
simple pole at $z'=0 $.
Thus, we arrive at the (dual vector space of the) space of global sections of $T^*C\otimes f^*l_j^* $ and of $T^* \Sigma'\otimes f^*l_j^* $ {\em with allowed simple poles}   at $w=0 $ and at $z'=0 $, respectively. These spaces of sections may be compared to the spaces $H^0(C, T^*C\otimes f^*l_j^*)^* $ and $H^0(\Sigma', T^* \Sigma'\otimes f^*l_j^*)^* $ via the exact sequences defined by the residue maps about $z'=0 $ and about $w=0 $, respectively.
Thus, we may consider (ii) as a result (see \cite{Brown}   for example) rather than as a definition. 

 {\bf 5.2.  A key ingredient of Theorem 2.}    Let $C' $ have the same meaning as in section 5.1 case i), and reserve the notation $C $ for case ii) except that the first marked point will also be allowed the role of nodal point of $\Sigma $ in $C' $.  The connected component of $[f;\Sigma, p(\Sigma)   ]$ in the space of $T $-fixed stable maps into $E_{\a(A) } $ is described as a fiber product of stable maps into the $T $-fixed strata of $E_{\a(A) } $.    A tree with root $C $ may connect, via a nodal point, to stable maps $C'\to\eps $ carrying the first marked point of $\Sigma $. 
 The smoothing of such a node  deforms $[f;\Sigma, p(\Sigma)   ]$ away from the locus of $T $-fixed stable maps into $E_{\a(A) } $. 
The inverse $T $-equivariant Euler class of the latter smoothing mode is given by $1/(-\psi_{\bullet}  +\chi_{\eps, \eps'}  /k)  $ where $\bullet $ is the smooth point of $C' $ in the normalization of $\Sigma $ that corresponds to the latter nodal point of $(f;\Sigma, p(\Sigma)  )  $.  Its presence is required by
the fixed-point localisation technique.  The tree with root $C $ yields a cohomology class of $B $ that is proportional to $1/(-z+\chi_{\eps, \eps'}  /k)  $ in contribution to the terms of type ii) in $\eps^*\J(z)  $. 
Let us observe that if we substitute $z\to \psi_{\bullet} $, then we get the $T $-equivariant Euler class $1/(-\psi_{\bullet}  +\chi_{\eps, \eps'}  /k)  $ of the latter smoothing mode. Let us integrate last over the moduli of $[f|_{C'}  , C', p(C')   ]$ where $C' $ is defined as in i).  The precedingly described nodal attachments to $C' $ yield terms of type  i) in ${\eps}  ^*\J(z) $.

If the tree with root $C $ is rooted at $\tilde \a $ there are two possible ways $f(C) $ can intersect with the stratum $\tilde \a $ at $f(0) $, according to the decomposition $\tilde \a^*= [1, \vec 0]^*+\tilde \pi^*\a^*pr_E $.  Namely, the pullback $[1, \vec 0]^* $ constrains $f(0) $ to lie in  $[1, \vec 0](A) $, while ${\tilde \pi}  ^* \a^*pr_E $ may be interpreted as constraining $f(0) $ to lie in ${\tilde \a}   \setminus  \P( N_{\a(A)}  \a(B)) $.
 
Define  \begin{align*}  t^{\eps}  (z)   :=\eps^*t(z)   -z+``\text{the sum of all}  &\ \text{contributions to ii)    where}  \\  \text{the first}\  & \text{marked point of $\Sigma $ is contained in $C $}  ".\end{align*}   
Let $\H^{\eps}  _+ $ be the completion of $\H^{(Euler_T^{-1}  (\cdot)   , \NN^{\eps}  ) }  _+ $ (Example 1.5)     by allowing additional additive terms that are infinite $z $ series at each order in Novikov's variables, of the form  $a\sum_{n=0}  ^{\infty}  (\chi/k)   ^{-b_n}  z^n $ where $b_n\geq n+1 $ and $a\in H^*(B, \QQ)  $.  Denote by $\J^{\eps} $ restrictions $\eps^*\J $ of $\J $ where $t^{\eps}  (z) $ is expanded in non-negative powers of $z $.
The preceding observations establish that $\J^{\eps}  (z, \tau)   \subset  \H^{\eps} $ is the point of the $(Euler_T^{-1}  (\cdot)   , \NN^{\eps}  )  $-twisted Lagrangian cone of $\eps $ with input $z+t^{\eps}  (z)   \in \H_+^{\eps} $.  Let us denote this Lagrangian cone contained in $\H^{\eps} $ by $\LL^{\eps} $.

{\bf 5.3. Recursion.}  

 The sheaf cohomology groups $H^i(C, f^*\NN^{\eps}  )  $ decompose as \begin{align*}  H^i(C, f^*\NN^{\eps}  )   =H^i(C, \oplus _j  f^*l_j)    \simeq&\\ \oplus _j H^i(C, f^*l_j)     \simeq& \oplus _j H^i(C, {\cal O}  ([f^*c_1(l_j)   ]\cap[C])  )   .\end{align*}  If $f^*c_1(l_j)   [
C]\geq 0 $, then a basis of $H^0(C, f^*l_j[C])  $ is given by $z^{f^*c_1(l_j)   [C]}  , \dots, 1 $.  In the present deformation problem $z $ is a section of the ${\cal O}  _C(1)  $ bundle, so $T $ acts on $z $ in the representation $\varphi^* $ with  weight $-\frac{\chi_{\eps, \eps'}}  {k} $.  
 The sections $z^{f^*c_1(l_j)   [C]}  , \dots, 1 $ may be realized explicitly as deformations of $[f; \Sigma , p(\Sigma)   ]$, as in \cite{Brown}.  Namely, each point of $E_{\a(A) } $ lies in either the normal bundle to the exceptional divisor, or its complement---both of which are toric bundles.

If $U_j(d_{\a, \b}  ) $ is negative, then we describe the virtual summand
 $H^1(C, f^* l_j) $ by the Serre duality theorem, $H^1 (C, f^*
l_j)\simeq H^0(C, T^*C \otimes f^*l_j^*)^* $.  If
\linebreak $(k, U_j(d_{\a, \b}  ))\neq (1, -1) $, then this is not the zero
vector space.  Otherwise, it is the $0 $-vector space, in which case it
does not contribute to $H^1(C, f^*TX) $.
 If $U_j(d_{\a, \b}  ) $ is negative, then
 $\{z^mw^{-kU(d_{\a, \b}  )-m}  \otimes\omega_C\}    _{m=1}  ^{-kU_j
(d_{\a, \b}  )-1} $ is a basis of the vector space $H^0 (C, T^*C \otimes
f^* l_j^*)\simeq H^1 (C, f^* l_j)^* $.  This set, when restricted to a
contractible neighborhood of $0\in C $, has the same weights as does
the basis $\{z^m \frac{dz}  {z}  \otimes\a^*l_j^*\}    _{m=1}  ^{-kU_j
(d_{\a, \b}  )-1} $ of the vector space $0^*H^0(C, T^*C\otimes
f^*l_j^*)\simeq 0^*H^1(C, f^*l_j)^* $.
  Finally, apply the Proposition describing the $l_j $ to compute
\begin{align*}  Coeff_{\eps, \eps'}  (k)=e^{\eps}  Euler_T ^{-1}     ({\mathcal N}  _{\eps, \eps'}  (k)   \oplus  
\NN^{\eps'}  )   .\end{align*}  
Given two of the $T $-fixed strata $\eps $ and $\eps' $ connected by an edge, define submanifolds of each where the strata intersect with edges connecting the two strata. The two submanifolds are diffeomorphic, call it $ Z_{\eps, \eps'} $, by the connecting edges. 
Let us point out the role of each of the terms of
\begin{align*}   \J^{\eps'}   (-\frac{\chi_{\eps, \eps'}}  {k}  , t)&=
-\frac{\chi_{\eps, \eps'}}  {k}  +\b^*t(\frac{\chi_{\eps, \eps'}}  {k}  )+\\&\eps'^*\sum_{n', D', d'}  \frac{Q^{D'}  q^{d'}}  {n'!}  (\ev_1)_*
\left[ \frac{1}  {-\frac{\chi_{\eps, \eps'}}  {k}  -\psi_1}  
\prod_{i=2}  ^{n'+1}   (\ev_i^*t)(\psi_i)\right]\end{align*}  on the RHS of the {\em recursion relation}  :
\begin{align*}   {\cal O}_{\eps, \eps'} \operatorname{Res}  _{z=-\frac{\chi_{\eps, \eps'}}  {k}}  \J^{\eps}  (z)\ dkz =&\\{\cal O}'_{\eps, \eps'} 
{\tilde q}  ^{k{\tilde d}  {\eps, \eps'}}   q^{kd_{\eps, \eps'}}  &Coeff_{\eps, \eps'}  (k)\
\J^{\eps'}  (-\frac{\chi_{\eps, \eps'}}  {k} ) , \end{align*} where 
\[{\cal O}_{\eps, \eps'}   =\left\{\begin{array}  {cl} Id  & \eps\simeq  \eps'\\  
\pi^{\perp} &\eps \simeq A, \eps'\simeq B \\  
  i_{ Z_{\eps, \eps'}}^* &\eps' \simeq A, \eps \simeq B \end{array}\right. \text{and}\  {\cal O}'_{\eps, \eps'}   =\left\{\begin{array}  {cl} Id  & \eps\simeq  \eps'\\  
i_{ Z_{\eps, \eps'}}^*  &\eps \simeq A, \eps'\simeq B \\  
 \pi^{\perp} &\eps' \simeq A, \eps \simeq B . \end{array}  \right.       \]

In order to proceed, choose a basis $\{\phi_{\mu}  \}   $ of $H^*(B) $, 
and denote by $\{\phi^{\mu}  \}   $ the Poincar\'e-dual  basis.  Let $\{
\xi_{\mu'}  \}   $ be a basis of $H^*(A) $ and $\{
\xi^{\mu'}  \}   $ dual basis vectors, which respects the orthogonal\footnote{with respect to $(\cdot, \cdot)_A $}     direct sum decomposition $H^*(A)\simeq Im(i_A^* )\oplus (Im(i_A^*))^{\perp} $. 
both direct summands are $Im(i_A^*) $--modules.  The role of ${\cal O}_{\eps, \eps'} $, ${\cal O}'_{\eps, \eps'} $ is understood by observing that 
 $Coeff_{(\a, j_+), \b}(k) $ are valued in $Im(i_A^*) $.   Finally, denote \[\phi_{(\mu, \mu')}  =\left\{\begin{array}  {cl}  \phi_{\mu}  & \eps\simeq B\\  
\xi_{\mu'}   &\eps\simeq A \end{array}  \right.        \text{and}  \    
\phi^{(\mu, \mu')}  =\left\{\begin{array}  {cl}  \phi^{\mu}  & \eps\simeq B\\  
\xi^{\mu'}   &\eps\simeq A. \end{array}  \right.       \]
 
Further, let us introduce ``delta-functions" at the $T $-fixed strata $\eps $: elements of $H_T^*(E_{\a(A)}  , \Q(\l)) $ given by

\[\delta_{\eps}  =\left\{\begin{array}  {ccl}   \frac{\prod_{j\notin\eps}  U_j}  {\prod_{j\notin\eps}  \eps^*U_j}   &&\text{if}  \ \eps  \simeq B\ \\   
&\frac{\prod_{jj\neq\eps}  U_{A, jj}}  {\prod_{jj\neq\eps}  \eps^*U_{A, jj}}  &\text{if}  \ \eps =[1, \vec 0]\\
\frac{-{\tilde P}  \prod_{jj\neq\eps}  U_{A, jj}}  {-\eps^*{\tilde P}  \prod_{jj\neq\eps}  \eps^*U_{A, jj}}  
&&\text{if}  \ \eps  \simeq A, \neq[1, \vec 0]\end{array}  \right.       \]
that satisfy \[\eps^*\delta_{\eps'}  =\delta_{\eps, \eps'}  1_{\eps}  \forall \eps, \eps'\in E_{\a(A)}  ^T.\]
It follows then, by fixed-point localization, that
 $\{\delta_{\eps}  \phi_{(\mu, \mu')}  \}   $ and $\{e^{\eps}  \delta_{\eps}  \phi^{(\mu, \mu')}  \}   $ are
Poincar\'e-dual bases of $H_T^*(E_{\a(A)}  , \Q(\l)) $.   Lastly, $\forall a\in
H^*((E_{\a(A)}  )_{0, n+1, \D}  ) $, 
\begin{align*}  (\ev_1)_*(a)=&\sum_{(\mu, \mu'), \eps} e^{\eps}  \delta_{\eps}  \phi^{(\mu, \mu')}  ((\ev_1)_*a, \delta_{\eps}  \phi_{(\mu, \mu')}  )_{E_{\a(A)}}  =\\&\sum_{(\mu, \mu'), \eps} e^{\eps}  \delta_{\eps}  \phi^{(\mu, \mu')}  (a, (\ev_1)^*(\delta_{\eps}  \phi_{(\mu, \mu')}  ))_{(E_{\a(A)}  )_{0, n+1, \D}}  =\\&\sum_{(\mu, \mu'), \eps} e^{\eps}  \delta_{\eps}  \phi^{(\mu, \mu')}  \int_{[(E_{\a(A)}  )_{0, n+1, \D}  ]}  a(\ev_
1)^*(\delta_{\eps}  \phi_{(\mu, \mu')}  ).
\end{align*}
  
Use the Poincar\'e pairing on $H^*(E_{\a(A)}  ) $ to process
\[ {\cal O}_{\eps, \eps'} (\ev_1)_*\left[ \frac{1}  {z-\psi_1}  \prod_{i=2}  ^{n+1}  
(\ev_i^*t)(\psi_i)\right]\]
according to the above formula.  Using fixed-point
localization to compute the integrals, we arrive at
\begin{align*}  {\cal O}_{\eps, \eps'}  e^{\eps}  \phi^{(\mu, \mu')}  \sum_{(\mu, \mu')}  \int_{[ Z_{\eps, \eps'}  ]\otimes
[(E_{\a(A)}  )_{0, n+1, \D-kd_{\eps, \eps'}}  ]}  \frac{\phi_{(\mu, \mu')}}  {k(z+\frac{\chi_{\eps, \eps'}}  {k}  )}  Euler_T^{-1}  ({\cal
N}  _{\eps , \eps'}  (k))\times \\ \frac{(\pi\ev_{\infty}  \times \pi\ev_0)^*([\Delta_{ Z_{\eps, \eps'}  \times 
 Z_{\eps, \eps'}}  ])}  {-\frac{\chi_{\eps, \eps'}}  {k}  -\psi_0}  \prod_{i=2}  ^{n+1}  
(\ev_i^*t)(\psi_i), \end{align*}  where $[\Delta_{ Z_{\eps, \eps'}  \times Z_{\eps, \eps'}}  ]$ stands for the fundamental class of the
diagonal in $ Z_{\eps, \eps'}  \times Z_{\eps, \eps'} $, the points $\infty\in C $ and $0\in\Sigma' $
lie in the normalization of $\Sigma $ over the nodal point where $C $
meets $\Sigma' $ in $\Sigma $, and we have 
evaluated $-\psi_{\infty}  =-\frac{\chi_{\eps, \eps'}}  {k} $.
      Rewrite this integral as follows, 
\begin{align*}  {\cal O}_{\eps, \eps'} e^{\eps}  \phi^{(\mu, \mu')}  \sum_{(\mu, \mu'), (\nu, \nu')}  \int_{[ Z_{\eps, \eps'}  ]}  \frac{\phi_{(\mu, \mu')}}  {k(z+\frac{\chi_{\eps, \eps'}}  {k}  )}  Euler_T^{-1}  ({\cal
N}  _{\eps, \eps'}  (k)\oplus 
N^{\eps'}  )e^{\eps'} \phi^{(\nu, \nu')}\times \\  \int_{[(E_{\a(A)}  )_{0, n+1, \D-kd_{\eps, \eps'}}  ]}  \frac{\ev_0^*(\delta_{\eps'}  \phi_{(\nu, \nu')}  )}  {-\frac{\chi_{\eps, \eps'}}  {k}  -\psi_0}  \prod_{i=2}  ^{n+1}  
(\ev_i^*t)(\psi_i), \end{align*}  to conclude with the terms of $\J^{\eps'}  (-\frac{\chi_{\eps, \eps'}}  {k}  ) $ that involve the push-forward
\[(\ev_1)_*:H^*((E_{\a(A)}  )_{0, n'+1, \D'}  )\rightarrow H^*(E_{\a(A)}  )\]on the RHS of the
recursion relation. In particular, we see the appearance of 
\begin{align*}  Coeff_{\eps, \eps'}  (k)
=e^{\eps}  Euler_T ^{-1}  ({\mathcal N}  _{\eps, \eps'}  (k)\oplus  
\NN^{\eps'}  ).\end{align*}  

The virtual summand $\ominus T_{\eps'}  X $ of ${\cal N}  _{\eps, \eps'}  (k) $
results from describing the virtual bundle
\[ H^0(\Sigma, f^*TX)\ominus H^1(\Sigma, f^*TX)\ominus (\Lie
\Aut(\Sigma; p(\Sigma))\cap H^0(\Sigma, f^*TX))\]
in terms of related virtual bundles  defined on the curve that
normalizes the node of $\Sigma $ where $C $ meets $\Sigma' $.  Hence, if
there is a marked point, rather than a node, of $\Sigma $ at $\infty\in
C $, then we add $\oplus T_{\eps'}  X $ to the formula for ${\cal
N}  _{\eps, \eps'}  (k) $.  If there is a node, rather than a marked point, of
 $\Sigma $  at $0\in C $, then we add $\ominus T_{\eps}  X $ to the
formula for ${\cal N}  _{\eps, \eps'}  (k) $.

Finally consider the case when a leg leaves a component of height 0, whose only marked or nodal point lies in the component of height 0.  This matches with the dilaton shift term of $\J^{\eps'}  (-\chi_{\eps, \eps'}  /k)  $ on the RHS of the recursion relation.

{\bf 5.4. Theorem 2.}  

{\tt Theorem 2.}   Points $\J(z)  $ of the overruled Lagrangian cone of the $T $-equivariant  genus-0 Gromov--Witten theory of $E_{\a(A) } $ are characterized by the conditions:

(1.a):   $\J^{(\a, jj_+) }  (-z)   \in \LL_{A}  ^{(\a, jj_+) } $   

(1.b):  $\J^{\c}  (-z)    \in \LL_B^{\c} $

(2.bb)   :
\[ \Res_{z=-\frac{\c^*U_{j_+}}  {k}}  \J^{\c}  (z)   dkz=q^{kd_{\c, \c'}}  Coeff_{\c, \c'}  (k)   \J^{\c'}  (-\frac{\c^*U_{j_+}}  {k}  )   \]

\begin{align*}  Coeff_{\c, \c'}  (k)   =\\&    \prod_{m=1}  ^{k-1}  
(m\frac{\chi_{\c, \c'}}  {k}  )   \prod_{m=1}  ^{k}  
(-m\frac{\chi_{\c, \c'}}  {k}  )   
\prod_{j\notin \c\cup\c'}  \prod_{m=1}  ^{kU_j(d_{\c, \c'}  ) }  
(\c^*U_j-m\frac{\chi_{\c, \c'}}  {k}  )   .\end{align*}  
 
(2.ab)   :

\[\pi^{\perp} \Res_{z=\frac{(\a, j_+)   ^*{\tilde P}}  {k}}  \J^{(\a, j_+) }  (z)   dkz=i_{A}  ^*{\tilde q}  ^{-k}  q^{kd_{\a, \b}}  Coeff_{(\a, j_+)   , \b}  (k)   
\J^{\b}  (\frac{(\a, j_+)   ^*{\tilde P}}  {k}  )   \]

\[   i_{A}  ^*\Res_{z=-\frac{\b^*U_{j_-}}  {k}} \J^{\b}  (z)   dkz=\pi^{\perp}{\tilde q}  ^{-k}  q^{kd_{\a, \b}}  Coeff_{\b, (\a, j_+) }  (k)   \J^{(\a, j_+) }  (-\frac{\b^*U_{j_-}}  {k}  )   
\]

\begin{align*}  Coeff_{(\a, j_+)   , \b}  (k)   
=&\\  i_A^*\prod_{m=1}  ^{k-1}  
(m\frac{\chi_{\a, \b}}  {k}  )   \prod_{m=1}  ^{k}  (-m\frac{\chi_{\a, \b}}  {k}  )   
&\prod_{j\notin \a\cup\b}  \prod_{m=1}  ^{kU_j(d_{\a, \b}  )   -k}  
({\tilde P}  ^{(\a, j_+) }  +\a^*U_j-m\frac{\chi_{\a, \b}}  {k}  )  \times \\
&\prod_{m=1}  ^k
(-{\tilde P}  ^{(\a, j_+) }  -m\frac{\chi_{\a, \b}}  {k}  )   
.\end{align*}

(2.aa)   :
 
\begin{align*}  \Res_{z=-\frac{(\a, jj_+)   ^*U_{A, jj_+'}}  {k}}  &\J^{(\a, jj_+) }  (z)   dkz=\\ {\tilde q}  ^{k}  &Coeff_{(\a, jj_+)   , (\a, jj'_+) }  (k)   \J^{(\a, jj'_+) }  (-\frac{(\a, jj_+)   ^*U_{A, jj_+'}}  {k}  )   
\end{align*}

\begin{align*}  Coeff_{(\a, jj_+)   , (\a, jj_+') }  (k)   =&\\  i_A^*\prod_{m=1}  ^{k-1}  
(m\frac{-\a^*U_{jj_+}  +\a^*U_{jj_+'}}  {k}  )  &\prod_{m=1}  ^k(-m\frac{-\a^*U_{jj_+}  +\a^*U_{jj_+'}}  {k}  )\times   \\
\prod_{jj\notin \a, jj\neq j  j_+'}  \prod_{m=1}  ^k
&({\tilde P}  ^{(\a, jj_+  )} +\a^*U_{jj}  -m\frac{-\a^*U_{jj_+}  +\a^*U_{jj_+'}}  {k}  )\times \\   
&\prod_{m=1}  ^{-k}   (-{\tilde P}  ^{(\a, jj_+) }  -m\frac{-\a^*U_{jj_+}  +\a^*U_{jj_+'}}  {k}  )   
\end{align*}

We have established that every point on $ \LL_{E_{\a(A) }} $ satisfies (1.a), (1.b) and (2).
Let us prove now that if series $\{ \J^{\eps}  \}   $ satisfy these conditions, 
then they represent a point in $\LL_{E_{\a(A) }} $. 

Denote by $-z+\eps^* t(z)  $ those terms in $t^{\eps}  (z)  $ that do not contribute to the residue (ii).
We will show that conditions (i)    and (ii)    for $\{\J^{\eps}  \}   $ describes an infinite sum over finite trees, with weights assigned to edges and vertices according to the same rules that one finds in fixed point localization of the point $\J(-z)  $ on $\LL_{E_{\a(A) }} $ whose projection to ${\mathcal H}  _+ $ is $\{-z+ \eps^* t(z)    \}   $.  Each finite tree will have a prescribed height, and will correspond to a particular element of $\H^{\eps} $.  Let ${\o \J}  ^{\eps}  (-z)  $ denote the point of $\LL^{\eps} $ that projects to $\eps^*t(z)   -z $ along $\H^{\eps}  _- $. 

Define the symbol $[\J^{\eps}  (-z)    ]^n $ inductively to be the sum over all weighted trees of height $n $, as follows.  Vertices are indexed by $T $-fixed $\eps''\in E_{\a(A) }  ^T $.  A vertex at height 1 indexed by $\eps'\in E_{\a(A) }  ^T $ can connect by an edge $e $ to a vertex of height 0 indexed by $\eps\in E_{\a(A) }  ^T $ provided that they are connected by an edge.  Associate an arbitrary positive integer $k_e\geq 1 $ to such an edge $e $ of the tree.  Associate the weight $q^{kd_{\eps, \eps'}}  Coeff_{\eps, \eps'}  (k_e)  $ to such an edge $e $, and the weight ${\o\J}  ^{\eps'}  (-\frac{\chi_{\eps, \eps'}}  {k}  )  $ to the vertex at height 1 indexed by $\eps'\in E_{\a(A) }  ^T $.  When there is only a vertex at height 0 at $\eps $, assign the weight ${\o \J}  ^{\eps}  (-z)   =:[\J^{\eps}  (-z)   ]^0 $.

When there is an edge $e $ connecting a vertex at height 1 at $\eps' $ to a vertex at height 0 at $\eps $, associate the weight $\frac{1}  {k_e(-z+\chi_{\eps, \eps'}  /k_e) } $ (expanded in positive powers of $z $)    to the vertex at height 0 at $\eps $.  Then combine the weights at height 0, the weight of the edge, and the weight at height 1 multiplicatively, and denote the result by $\f^{k_e}  _{\eps, \eps'}  (-z)  $.  This assignment of weights is uniquely implied by (ii), when we truncate $\J^{\eps'}  (-\chi_{\eps, \eps'}  /k)  $ (within $\LL^{\eps} $)    on the RHS by replacing $t^{\eps'}  (-\chi_{\eps, \eps'}  /k)  $ by $\eps'^*t(-\chi_{\eps, \eps'}  /k)  $.  To see this, let  us point out that $-z+\eps^*t(z)  $ is the truncation of $t^{\eps}  (z)  $ modulo  its summands that contribute to the residues ii). 

Let ${\o \J}  ^{\eps, 1}  (-z)  $ denote the point of $\LL^{\eps} $ that projects to $\sum_{\eps'}  \sum_{k_e>0}  \f^{k_e}  _{\eps, \eps'}  (-z)+\eps^*t(z)   -z $ along $\H^{\eps}  _- $, and define $[\J^{\eps}  (-z)   ]^1:={\o \J}  ^{\eps, 1}  (-z)   -\eps^*t(z)+z $.  The terms $\eps^*t(z)   -z $ are subtracted to prevent duplication with $[\J^{\eps}  (-z)   ]^0 $.  Suppose $[\J^{\eps}  (-z)    ]^{n-1} $ has already been defined and has been given an interpretation as a sum over weighted trees, with the vertices of height $n-1 $ at $\eps'' $ (connected by an edge $e $ to vertices at height $n-2 $ at $\eps''' $)    assigned to weights ${\o \J}  ^{\eps''}  (-\chi_{\eps''', \eps''}  /k_e)  $.  Then define $[\J^{\eps}  (-z)    ]^n $ by using the same weighted trees as for $[\J^{\eps}  (-z)    ]^{n-1} $, but with the vertices of height $n-1 $ at $\eps'' $ assigned to the {\em different}   weights ${\o \J}  ^{\eps'', 1}  (-\chi_{\eps''', \eps''}  /k_e)  $.  
 For each $n\geq 0 $, the vectors $[\J^{\eps}  (-z)    ]^n $ by their very construction are contained in $\LL^{\eps} $, satisfy the residue condition (ii)    up to height $n-1 $, and are uniquely determined by $\{\eps''^*t(z)   \}    _{\eps''\in E_{\a(A) }  ^T} $.

Let us now turn to fixed point localization of the point $\J(-z)  $ on $\LL_{E_{\a(A) }} $ whose projection to ${\mathcal H}  _+ $ is $\{-z+\eps^* t(z)    \}   $.  Consider any $T $-fixed stable map $(f; \Sigma)  $ that contributes to $\eps^*\J (-z)  $. Consider any component $C''\subset \Sigma $ that maps to a $T $-fixed section of $E_{\a(A) } $. We define the height $n $ of $C'' $ to be the unique number of distinct legs of $(f; \Sigma)  $ that must be traversed by any path starting at  the first marked point of $\Sigma $, andending in $C'' $.  We define the height $n $ of a leg to be the unique number of legs that must be intersected by a path starting at the first marked point of $\Sigma $ andending in the interior of the leg.

We will consider an increasing filtration on the set of $T $-fixed stable maps that contribute to $\eps^*\J (-z)  $.  Define the $n^{th} $ filtered part, $[\eps^*\J (-z)   ]^{n} $, to consist of all stable maps whose  components  have height at most $n $.   

It follows immediately that for each $n\geq 0 $, the following elements of $H^{\eps} $ coincide, for a given initial condition $\eps^* t(z)  $---the contribution to fixed point localization of Gromov--Witten invariants in $\eps^*\J (-z)  $ due to elements of $[\eps^*\J (-z)   ]^n $, and the sum of weighted trees constructed by (i)    and (ii)    that contribute to $[\J^{\eps}  (-z)    ]^n $.

\section{Recursion}

To prove the equivariant version of Theorem, it suffices to show that 
 $\J=I_{E_{\a(A) }} $ satisfies conditions (1.a), (1.b) and (2)    of Theorem 2.
Define \[{\mathscr U}  _{\cal J}  =\left\{\begin{array}  {lll}  U_j& , {\cal J}   =j\\U_{A, jj}  &, {\cal J}   =A, jj\\{\tilde P}  &, {\cal J}  =\text{``else"}  \end{array}  \right.        \in H^2_T(E_{\a(A) }  )   .\]
The hypergeometric modification $I_{E_{\a(A) }} $ is a $(q, {\tilde q}  , Q)  $-series whose coefficients have simple poles at $z=-\eps^*{\mathscr U}  _{\cal J}  
/k $, finite order poles at $z=\infty $, and essential
singularities at $z=0 $. Thus, we need to show
that: (1.a)
 $(\a, jj_+)^*I_{E_{\a(A) }}  \in \LL^{(\a, jj_+)} $, 
 (1.b) $\c^*I_{E_{\a(A) }}  \in \LL^{\c} $, and (2)    residues at the simple poles
satisfy the recursion relation of Theorem 2. We postpone (1.a), (1.b) until Section 7, 
and deal with (2)    here by computing the residues explicitly.

Our first goal is to argue that the series $\eps^*I_{E_{\a(A) }} $ is  supported in the Mori cone of $E_{\a(A) }  , \forall\eps \neq{\tilde \a}  (B\setminus A) $.  The mechanism that insures this is to look at the support of the factors $\frac{1}   {\prod_{m=1}  ^{{\mathscr U}  _{\mathcal J}  ({\tilde {\cal D}}  ) }  ({\mathscr U}  _{\mathcal J}  +mz) } $ of $I_{E_{\a(A)}}  (z) $ for which $\eps^*{\mathscr U}  _{\mathcal J}  =0 $, $\forall \mathcal J\in\eps, \forall\eps \neq{\tilde \a}  (B\setminus A) $.

{\tt Proposition.}   {\em Any element of $MC(E_{\a(A)}  )  $ may be represented by a curve whose irreducible components are preserved by the action of $T $ on $E_{\a(A)} $.}  

{\em Proof.}    The action of $T $ on $E_{\a(A)} $ is induced by that on $E $.  Thus, we need only prove the result for $E $. Extend the rows of $(m_{ij}  )  $ to an integer basis of ${\Z}  ^N $, each defining a unitary circle action $S_l^1:E\rightarrow E $ and a complexified circle action $\C^*_l:E\rightarrow E $.  For each $p\in E $, denote $(X_{l, \theta}  )   _p $ the tangent vector at $p $ induced by the $S_l^1 $-action, and $(X_{l, r}  )   _p:= -J_p (X_{l, \theta}  )   _p $.  Recall that $\omega_p(\cdot, (X_{l, \theta}  )   _p)   =d\mu_p (X_{l, \theta}  )   (\cdot)  $.
Therefore, $\mu_p((X_{l, \theta}  )   _p)  $ is an increasing function along the flow lines of $X_{l, r} $ and is constant along the flow lines of the remaining $N-K-1 $ radial vector fields.  Since the fibers of $E $ are compact, it follows that the flow lines of $X_{l, r} $ must come to an end asymptotically, and this happens at the points where $X_{l, r} $ vanishes. 

Gromov's compactness theorem \cite{Gromov}   guarantees that the limit set of the flow of such a curve is a (possibly nodal)    curve in $E $.  Iterating this argument for each of the flows, the result follows.

With the Proposition in place, let us now compare $MC(E_{\a(A)}  ) $ and \linebreak
 $MC(E_{\a(B)}  ) $.  A first source of difference between the two comes from the inclusion
 $(\a, jj_+)_*{i_A}  _* MC(A)\subset  (\a, jj_+)_*MC(B) $.  
Another difference is that the $T $-invariant $\C P^1_{\a(B\setminus A), \b(B\setminus A)} $ curves in $E_{\a(A)} $ do not have any geometric analogues in $E_{\a(B)} $.       However, the latter curve may be represented as the sum of the class of a $\C P^1_{(\a, j_+), \b} $ and the class of a $\C P^1 $ in a fiber of the exceptional divisor.  Thus all elements of $Ker\pi_* $ have geometric analogues in $H_2(\text{fiber of}  \ E_{\a(B)}  ) $.

{\em  Remark.}     Any curve from a fiber of $E_{\a(B)} $ has a geometric analogue in $\pi^{-1}  (A) $.  The $U_j $'s are determined by the geometry of $E $, and thus have the same meaning whether pulled back to $H^2_T (\text{fiber of}  \ \pi^{-1}  (A)) $ or to \linebreak
 $H^2_T(\text{fiber of}  \ E_{\a(B)}  ) $.
 The class $\tilde P $ is determined\footnote{As a functional on the classes of $T $-invariant curves.}   by the local geometry of the exceptional divisor and thus has the same meaning whether referred  to $H^2_T(\text{fiber of}  \ \pi^{-1}  (A)) $ or to $H^2_T(\text{fiber of}  \ E_{\a(B)}  ) $.

\noindent{\bf  (2.aa)    Residue of $\J^{(\a, jj_+) }  (z)  $ at $z=-\frac{\chi_{(\a, jj_+)   , (\a, jj'_+) }}  {k}  , \ k\geq 1 $.}   Given $D\in MC(B)   , $
rename $d_i\to d_i' $ and ${\tilde d}  \to {\tilde d'} $, and then redefine $d_i' $ and ${\tilde d'} $,

\[{\tilde d}  '={\tilde P}  ^{(\a , jj_+) }  (D)+{\tilde d}  , \]
  
\[ d'=P^{\a}  (D)+d.\]
Then, the pairings $U_j({\cal D}  )   =U_j(d)   -\L_j(D)  $ translate into $U_j(d)+\a^*U_j(D)  $.
\begin{align*}  (\a, jj_+)   ^*I_{E_{\a(A) }}  =i_A^*e^{{\tilde P}  ^{(\a, jj_+) }  {\tilde t}  /z} e^{P^{\a}  t/z}  \sum_{D\in MC(B) }  
\sum_{d \in \Z^K, {\tilde d}  \in\Z}   \\
\frac{({\bf z^{ht_{A  }  (D') }}  \J^D_{e(\cdot), L}  (z, \tau) +\G^D_{ A}  (z, \tau))  (Q^D{\tilde q}  ^{{\tilde P}  ^{(\a, jj_+)}     (D) }  q^{P^{\a}  (D) }  )q^d{\tilde q}  ^{{\tilde d}} e^{{\tilde d}  {\tilde t}} e^{dt} e^{{\tilde P}  ^{(\a, jj_+)}     (D)   {\tilde t}} e^{P^{\a}  (D)   t}}  
{\prod_{j\notin \{j_+(\a, \b)   |  \a\sim\b\}  }  
\prod_{m=1}  ^{U_j(d)+\a^*U_j(D) }  (\a^*U_j+mz) }  \times \\
\frac{1}  
{ \prod_{j\in \{j_+(\a, \b)   |  \a\sim\b\}  }  
\prod_{m=1}  ^{U_j(d)+\a^*U_j(D)+{\tilde P}  ^{(\a, jj_+)}    (D)+{\tilde d}}  (\a^*U_j+{\tilde P}  ^{(\a, jj_+) }  +mz) }  \times \\
\frac{1}  {\prod_{m=1}  ^{-{\tilde d}  -{\tilde P}  ^{(\a, jj_+)}   (D) }  (-{\tilde P}  ^{(\a, jj_+) }  +mz)   
}  \times \\ \frac{1}  {\prod_{m=1}  ^{c_1(L)   (D)+{\tilde P}  ^{(\a, jj_+)}   (D)+{\tilde d}}  (c_1(L)+{\tilde P}  ^{(\a, jj_+) }  +mz) }  .\end{align*}  
{\tt Proposition.}  {\em The series $(\a, jj_+)   ^*I_{E_{\a(A) }} $ is  supported in the Mori cone of $E_{\a(A) } $.}   

{\em Proof.}    For $j_+:=jj_+\neq [1, \vec 0]$, the support of the series $\frac{1}  {\prod_{m=1}  ^{U_{j_+(\a, \b)}  (d)+\tilde d}   mz} $ is characterised by the inequality ${U_{j_+(\a, \b)}  (d)+\tilde d}  \geq 0 $.  For each $j\in\a $, the support of the series $\frac{1}  {\prod_{m=1}  ^{U_j(d) }  mz} $ is characterised by the inequality $U_j(d)   \geq 0 $.  Let us now argue that the set of solutions $(d, \tilde d) $ to the same inequalities is contained in $Ker \pi_* $. 
By the comparison of $Ker\pi_* $ with $H_2(\text{fiber of}  \ E_{\a(B)}  ) $, and by the Remark, it suffices to establish the analogous result for 
 $H_2(\text{fiber of}  \ E_{\a(B)}  ) $.  This follows from the Corollary and the same (strictly speaking, analogous) inequalities that arise there, as a special case of a general result in toric geometry describing the Mori cone in terms of inequalities. 

For $jj_+=[1, \vec 0]$, the inequalities describing the support of the series are $U_j(d)\geq 0\forall j\in\a $ and ${\tilde d}  \geq 0 $, whose solution set is ``a subset of $MC(\text{fiber of}  \ E)\linebreak\subset  MC(E_{\a(A)}  ) $" $\oplus \N\cdot $``the class of a $\C P^1 $ in a fiber of the exceptional divisor".

Let us now factor and rewrite the terms of the residue of the above pullback series at $z=-\chi/k $ as follows:
\begin{align*}   &{\bf If\ jj_+=[1, \vec 0]:}    \prod_{m=1}  ^{c_1(L)   (D)+{\tilde P}  ^{(\a, jj_+)}    (D)+{\tilde d}}  (c_1(L)+{\tilde P}  ^{(\a, jj_+) }  +m(-\chi/k)  )   =\\& \prod_{m=1}  ^{k}  ( 0   -m\frac{\chi}  {k}  )   \times \\&
\prod_{m=1}  ^{c_1(L)   (D)+{\tilde P}  ^{(\a, jj'_+) }  (D)+{\tilde P}  ^{(\a, jj_+) }  (D)   -{\tilde P}  ^{(\a, jj'_+) }  (D)+{\tilde d}  -k }  (c_1(L)+{\tilde P}  ^{(\a, jj'_+) }  -m\frac{\chi}  {k}  ), \\
& \text{for}   \ jj'_{+}  =  j\in \{j_+(\a, \b)   |  \a\sim\b\}    , 
\prod_{m=1, \neq k}  ^{U_j(d)+\a^*U_j(D)+{\tilde P}  ^{(\a, jj_+)}     (D)+{\tilde d}}  (\a^*U_j+{\tilde P}  ^{(\a, jj_+) }  -m\frac{\chi}  {k}  )   =
 \\&
\prod_{m=1}  ^{k-1}  
 (\chi-m\frac{\chi}  {k}  )   \times \\ 
& \prod_{m=1}  ^{U_j(d)+\a^*U_j(D)+{\tilde P}  ^{(\a, jj'_+) }  (D)+{\tilde P}  ^{(\a, jj_+) }  (D)   -{\tilde P}  ^{(\a, jj'_+) }  (D)   
+{\tilde d}  -k }  (\a^*U_j+{\tilde P}  ^{(\a, jj'_+) }  -m\frac{\chi}  {k}  ).\\&{\bf If\ jj_+\neq [1, \vec 0]:}  \end{align*}   \begin{align*}  &{\bf If\ jj'_+=[1, \vec 0]}       , \prod_{m=1, \neq  k}  ^{c_1(L)   (D)+{\tilde P}  ^{(\a, jj_+)}   (D)+{\tilde d}}  (c_1(L)+{\tilde P}  ^{(\a, jj_+) }  +m(-\chi/k)  )   = \\& \prod_{m=1}  ^{k-1}  (\chi-m\frac{\chi}  {k}  )   \times \\&
\prod_{m=1}  ^{c_1(L)   (D)+{\tilde P}  ^{(\a, jj'_+) }  (D)+{\tilde P}  ^{(\a, jj_+) }  (D)   -{\tilde P}  ^{(\a, jj'_+) }  (D)+{\tilde d}  -k }  (c_1(L)+{\tilde P}  ^{(\a, jj'_+) }  -m\frac{\chi}  {k}  ), 
\\
&\prod_{m=1}  ^{U_{jj_+}  (d)+\a^*U_{jj_+}  (D)+{\tilde P}  ^{(\a, jj_+)}     (D)+{\tilde d}}  (\a^*U_{jj_+}  +{\tilde P}  ^{(\a, jj_+) }  -m\frac{\chi}  {k}  )   = \prod_{m=1}  ^{k}  
 (\chi-m\frac{\chi}  {k}  )   \times \\
&\prod_{m=1}  ^{U_{jj_+}  (d)+\a^*U_{jj_+}  (D)+{\tilde P}  ^{(\a, jj'_+) }  (D)+{\tilde P}  ^{(\a, jj_+) }  (D)   -{\tilde P}  ^{(\a, jj'_+) }  (D)   
+{\tilde d}  -k }  (\a^*U_{jj_+}  +{\tilde P}  ^{(\a, jj'_+) }  -m\frac{\chi}  {k}  ). \\&{\bf If\ jj'_+=j'(\a, jj_+)\ (\neq [1, \vec 0])}  , \ \ \prod_{m=1, \neq k}  ^{U_{j'}  (d)+\a^*U_{j'}  (D)+{\tilde P}  ^{(\a, jj_+)}    (D)+{\tilde d}}  (\a^*U_{j'}  +{\tilde P}  ^{(\a, jj_+) }  -m\frac{\chi}  {k}  )   =  \\&\prod_{m=1}  ^{k-1}  
 (\chi-m\frac{\chi}  {k}  )   \times \\
&\prod_{m=1}  ^{U_{j'}  (d)+\a^*U_{j'}  (D)+{\tilde P}  ^{(\a, jj'_+) }  (D)+{\tilde P}  ^{(\a, jj_+) }  (D)   -{\tilde P}  ^{(\a, jj'_+) }  (D)   
+{\tilde d}  -k }  (\a^*U_{j'}  +{\tilde P}  ^{(\a, jj'_+) }  -m\frac{\chi}  {k}  ), \\&     \prod_{m=1}  ^{U_{jj_+}  (d)+\a^*U_{jj_+}  (D)+{\tilde P}  ^{(\a, jj_+)}   (D)+{\tilde d}}  (\a^*U_{jj_+}  +{\tilde P}  ^{(\a, jj_+) }  -m\frac{\chi}  {k}  )   =   \prod_{m=1}  ^{k}  
 (\chi-m\frac{\chi}  {k}  )   \times \\
&\prod_{m=1}  ^{U_{jj_+}  (d)+\a^*U_{jj_+}  (D)+{\tilde P}  ^{(\a, jj'_+) }  (D)+{\tilde P}  ^{(\a, jj_+) }  (D)   -{\tilde P}  ^{(\a, jj'_+) }  (D)   
+{\tilde d}  -k }  (\a^*U_{jj_+}  +{\tilde P}  ^{(\a, jj'_+) }  -m\frac{\chi}  {k}  ). \end{align*}   The terms that come before the `` $\times $" sign are associated with the deformations along $T\C P^1_{\eps, \eps'} $.  The remaining factors\footnote{When a factor of the series $Res_{z=-\chi/k}  I_{E_{\a(A)}}  (z) dkz $ is repeated, as regards its appearance in the above set of cases of values of $jj_+ $ and $jj_+' $, use the transformation laws of the above cases.}   transform as:\begin{align*}  
 &  \prod_{m=1}  ^{c_1(L)   (D)+{\tilde P}  ^{(\a, jj_+)}    (D)+{\tilde d}}  (c_1(L)+{\tilde P}  ^{(\a, jj_+) }  +m(-\chi/k)  )   =\\ & \prod_{m=1}  ^{k}  (   c_1(L)+{\tilde P}  ^{(\a, jj_+)}    -m\frac{\chi}  {k}  )   \times   \\ &
\prod_{m=1}  ^{c_1(L)   (D)+{\tilde P}  ^{(\a, jj'_+) }  (D)+{\tilde P}  ^{(\a, jj_+) }  (D)   -{\tilde P}  ^{(\a, jj'_+) }  (D)+{\tilde d}  -k }  (c_1(L)+{\tilde P}  ^{(\a, jj'_+) }  -m\frac{\chi}  {k}  ), \\
&\text{For}   \ j\in \{j_+(\a, \b)   |  \a\sim\b\}    , \neq jj'_+, \\&\prod_{m=1}  ^{U_j(d)+\a^*U_j(D)+{\tilde P}  ^{(\a, jj_+)}    (D)+{\tilde d}}  (\a^*U_j+{\tilde P}  ^{(\a, jj_+) }  -m\frac{\chi}  {k}  )   =
 \prod_{m=1}  ^{k}  
 (\a^*U_j+{\tilde P}  ^{(\a, jj_+) }  -m\frac{\chi}  {k}  )   \times \\
& \prod_{m=1}  ^{U_j(d)+\a^*U_j(D)+{\tilde P}  ^{(\a, jj'_+) }  (D)+{\tilde P}  ^{(\a, jj_+) }  (D)   -{\tilde P}  ^{(\a, jj'_+) }  (D)   
+{\tilde d}  -k }  (\a^*U_j+{\tilde P}  ^{(\a, jj'_+) }  -m\frac{\chi}  {k}  ) , 
\end{align*}  \begin{align*}   
&\text{for}   \ j\notin \{j_+(\a, \b)   |  \a\sim\b\}    
, \ \ \prod_{m=1}  ^{U_j(d)+\a^*U_j(D)}  (\a^*U_j-m\frac{\chi}  {k}  )   =  1 \times \\
 &  \prod_{m=1}  ^{U_j(d)+\a^*U_j(D)}  (\a^*U_j-m\frac{\chi}  {k}  ), \\
&\prod_{m=1}  ^{-{\tilde P}  ^{(\a, jj_+) }  (D)   -{\tilde d}}  (-{\tilde P}  ^{(\a, jj_+) }  -m\frac{\chi}  {k}  )   =\prod_{m=1}  ^{k}  (-{\tilde P}  ^{(\a, jj_+) }  -m\frac{\chi}  {k}  )   \times \\
&\prod_{m=1}  ^{-{\tilde P}  ^{(\a, jj'_+) }  (D)   -{\tilde P}  ^{(\a, jj_+) }  (D)+{\tilde P}  ^{(\a, jj'_+) }  (D)   -{\tilde d}  +k }  (-{\tilde P}  ^{(\a, jj'_+) }  -m\frac{\chi}  {k}  ), 
\end{align*}  \begin{align*}  & (Q^Dq^{P^{\a}  (D) }  {\tilde q}  ^{{\tilde P}  ^{(\a, jj_+)}     (D) }  )    q^d {\tilde q}  ^{\tilde d}  = {\tilde q}  ^k \times (Q^D q^{P^{\a}  (D) }  {\tilde q}  ^{{\tilde P}  ^{(\a, jj'_+)}     (D) }  )   \
\\&q^d{\tilde q}  ^{{\tilde d}  -k+{\tilde P}  ^{(\a, jj_+)}     (D)   -{\tilde P}  ^{(\a, jj'_+)}     (D) }  , \\
&\exp   \left( -\frac{(P^{\a}  t+{\tilde P}  ^{(\a, jj_+) }  {\tilde t}  )   k}  {\chi}  \right)   
\exp (dt+{\tilde d}  {\tilde t}  )    \exp (P^{\a}  (D)   t+{\tilde P}  ^{(\a, jj_+) }  (D)   {\tilde t}  )   =  1  \times \\&
\exp \left(-\frac{(P^{\a}  t+{\tilde P}  ^{(\a, jj'_+) }  {\tilde t}  )   k}  {\chi}  \right)   
\exp \left[dt+\left({\tilde d}  -k+{\tilde P}  ^{(\a, jj_+) }  (D)   -{\tilde P}  ^{(\a, jj'_+) }  (D)   \right)   {\tilde t}  \right]
\\&\exp (P^{\a}  (D)   t+{\tilde P}  ^{(\a, jj'_+) }  (D)   {\tilde t}  )   .
\end{align*}  In the last
equality we use $({\tilde P}  ^{(\a, jj_+) }  -{\tilde P}  ^{(\a, jj'_+) }  )/\chi=1 $
 from 2.2. 

Factors on the R.H.S. which come before the multiplication
sign `` $\times $'' form the recursion coefficients
 ${\tilde q}  ^kCoeff_{(\a, jj_+)   , (\a, jj'_+) }  (k)  $. Factors which
come after the multiplication sign form the term of the series $\J^{(\a, jj'_+) } $
evaluated at $z=-\chi/k $ and with the summation index $\tilde d $ replaced
with ${\tilde d}  -k+{\tilde P}  ^{(\a, jj_+) }  (D)   -{\tilde P}  ^{(\a, jj'_+) }  (D)  $.   Reversing this change in the 
summation index, we conclude that 
\begin{align*}  \Res_{z=-\frac{\chi_{(\a, jj_+)   , (\a, jj'_+) }}  {k}}  \ \J^{(\a, jj_+) }  (z)   \ dkz =&\\ {\tilde q}  ^k
{Coeff_{(\a, jj_+)   , (\a, jj'_+) }  (k) }  \ \J^{(\a, jj'_+) }  &(-\frac{\chi_{(\a, jj_+)   , (\a, jj'_+) }}  {k}  )   , \end{align*}  as required.

\noindent{\bf  (2.bb)    Residue of $\J^{\c}  (z)  $ at $z=-\frac{\chi_{\c, \c'}}  {k}  , \ k\geq 1 $.}   Given $D\in MC(B)   , $ and $\c\neq \c' $, 
rename $d_i\to d_i' $ and ${\tilde d}  \to {\tilde d'} $, and then redefine $d_i' $ and ${\tilde d'} $,

\[{\tilde d}  '={\tilde P}  ^{\c}  (D)+{\tilde d}  , \]

\[ d'=P^{\c}  (D)+d.\]
The pullbacks $\c^*{\tilde P} $ vanish.  In particular ${\tilde q}  ^{{\tilde P}  ^{\c}  (D) }  =1 $, and

\[{\tilde d}  '={\tilde d}   , \]

\[ d'=P^{\c}  (D)+d.\]

Then, the pairings $U_j({\cal D}  )   =U_j(d)   -\L_j(D)  $ translate into $U_j(d)+\c^*U_j(D)  $.
\begin{align*}  \c^*I_{E_{\a(A) }}  =e^{P^{\c}  t/z}   \sum_{D\in MC(B) }  \sum_{d\in \Z^K
, {\tilde d}  \in\Z}  z^{ht_A(D')}     \times \\
\frac{\J_B^D(z, \tau)    (Q^Dq^{P^{\c}  (D) }  )    q^d{\tilde q}  ^{{\tilde d}} e^{{\tilde d}  {\tilde t}} e^{dt} e^{P^{\c}  (D)   t}}  
{ \prod_{j\in \{j_+(\a, \b)   |  \a\sim\b\}  }  \prod_{m=1}  ^{U_j(d)+\c^*U_j(D)+{\tilde d}}  (\c^*U_j+mz) }  \times \\
\frac{1}  {\prod_{j\notin \{j_+(\a, \b)   |  \a\sim\b\}  }  \prod_{m=1}  ^{U_j(d)+\c^*U_j(D) }  (\c^*U_j+mz) }  \times \\   \frac{1}  {\prod_{m=1}  ^{-{\tilde d}}  (mz)   \prod_{m=1}  ^{{\tilde d}}  (c_1(L)+c_1(L)   (D)   z+mz) }  .\end{align*}  {\tt Proposition.}  {\em If $\c\neq{\tilde \a}  (B\setminus A) $ then the series $\c^*I_{E_{\a(A) }} $ is  supported in the Mori cone of $E_{\a(A) } $.}  

{\em Proof.}   The support of the series $\frac{1}  {\prod_{m=1}  ^{-\tilde d}   mz} $ is characterised by the inequality ${-\tilde d}  \geq 0 $.  
The set $\{j_+(\a, \b)   | \b\sim \a\}   $ coincides with the set $\{j\notin\a \}   $.  For each $j\in\a^c\cap\c $, the support of the series $\frac{1}  {\prod_{m=1}  ^{U_j(d)+\tilde d}   mz} $ is characterised by the inequality $U_j(d )+{\tilde d}  \geq 0 $.    For each $j\in\a\cap\c $, the support of the series $\frac{1}  {\prod_{m=1}  ^{U_j(d) }  mz} $ is characterised by the inequality $U_j(d)   \geq 0 $. 
The proof proceeds as in the case of 2.aa $(jj_+\neq [1, \vec 0]) $.

         {\tt Remark.}   For $\c= {\tilde \a}  (B\setminus A) $, the inequalities describing the support of the series are $U_j(d)\geq 0\forall j\in\a $ and ${\tilde d}  \leq 0 $, whose solution set is ``a subset of $MC(\text{fiber of}  \ E)\subset MC(E_{\a(A)}  ) $" $\oplus -\N\cdot $``the class of a $\C P^1 $ in a fiber of the exceptional divisor". 

Since ${\tilde P}  ^{\c}  =0 $ for all $\c\in E^T $ it follows that ${\tilde P}  (d_{\c, \c'}  )   =0 $.  Hence the ``index" ${\tilde d} $ does not transform presently.  Thus, the asymmetry between the factors indexed by $ j\notin\{j_+(\a, \b)   | \b\sim \a\}   $ and $ j\in\{j_+(\a, \b)   | \b\sim \a\}   $ is removed for the purposes of the present recursion process.  It follows that the present recursion process is identical to the toric bundles case \cite{Brown}, 

\[ \Res_{z=-\frac{\chi_{\c, \c'}}  {k}}  \ \J^{\c}  (z)   \ dkz = q^{kd_{\c, \c'}}  
Coeff_{\c, \c'}  (k)   \ \J^{\c'}  (-\frac{\chi_{\c, \c'}}  {k}  )   , \]
as required.

\noindent{\bf  (2.ab)    Residue of $\J^{(\a, j_+(\a, \c)  ) }  (z)  $ at $z=-\frac{\chi_{(\a, j_+(\a, \c)  )   , \c}}  {k}  , \ k\geq 1 $.}   Given $D\in MC(B)   , $
rename $d_i\to d_i' $ and ${\tilde d}  \to {\tilde d'} $, and then redefine $d_i' $ and ${\tilde d'} $, 
\[{\tilde d}  '={\tilde P}  ^{(\a, j_+) }  (D)+{\tilde d}  , \]

\[ d'=P^{\a}  (D)+d.\]
Then, the pairings $U_j({\cal D}  )   =U_j(d)   -\L_j(D)  $ translate into $U_j(d)+\a^*U_j(D)  $.

\begin{align*}  (\a, j_+(\a, \c)  )   ^*I_{E_{\a(A) }}  =i_A^*e^{{\tilde P}  ^{(\a, j_+(\a, \c)  ) }  {\tilde t}  /z} e^{P^{\a}  t/z}  
\sum_{D\in MC(B) }  
\sum_{d \in \Z^K, {\tilde d}  \in\Z}    \\
\frac{({\bf z^{ht_{A  }  (D') }}  \J^D_{e(\cdot), L}  (z, \tau) +\G^D_{A}  (z, \tau))(Q^D{\tilde q}  ^{{\tilde P}  ^{(\a, j_+(\a, \c)  ) }  (D) }  q^{P^{\a}  (D) }  )    q^d{\tilde q}  ^{{\tilde d}} e^{{\tilde d}  {\tilde t}} e^{dt} e^{{\tilde P}  ^{(\a, j_+(\a, \c)  ) }  (D)   {\tilde t}} e^{P^{\a}  (D)   t}}  
{\prod_{j\notin \{j_+(\a, \b)   |  \a\sim\b\}  }  
\prod_{m=1}  ^{U_j(d)+\a^*U_j(D) }  (\a^*U_j+mz) }  \times \\
\frac{1}  
{ \prod_{j\in \{j_+(\a, \b)   |  \a\sim\b\}  }  
\prod_{m=1}  ^{U_j(d)+\a^*U_j(D)+{\tilde P}  ^{(\a, j_+(\a, \c)  )}     (D)+{\tilde d}}  (\a^*U_j+{\tilde P}  ^{(\a, j_+(\a, \c)  ) }  +mz) }  \times \\
\frac{1}  {\prod_{m=1}  ^{-{\tilde d}  -{\tilde P}  ^{(\a, j_+(\a, \c)  ) }  (D) }  (-{\tilde P}  ^{(\a, j_+(\a, \c)  ) }  +mz)   
}  \times \\  \frac{1}  {\prod_{m=1}  ^{c_1(L)   (D)+{\tilde P}  ^{(\a, j_+(\a, \c)  ) }  (D)+{\tilde d}}  (c_1(L)+{\tilde P}  ^{(\a, j_+(\a, \c)  ) }  +mz) }  .\end{align*}  
Let us now factor and rewrite the terms of the residue of the above pullback series at $z=-\chi/k $ as follows:

\begin{align*}   
&\prod_{m=1, \neq k}  ^{-{\tilde P}  ^{(\a, j_+(\a, \c)  ) }  (D)   -{\tilde d}}  (-{\tilde P}  ^{(\a, j_+(\a, \c)  ) }  -m\frac{\chi}  {k}  )   =
\prod_{m=1}  ^{k-1}  (\chi-m\frac{\chi}  {k}  )   \times \\
&\prod_{m=1}  ^{\c^*U_j(D)+\a^*U_j(D)   -\c^*U_j(D)   -{\tilde P}  ^{(\a, j_+(\a, \c)  ) }  (D)+0-{\tilde d}  +k }  (-0-m\frac{\chi}  {k}  )\\
&\prod_{m=1}  ^{U_{j_-(\a, \c)}  (d)+\a^*U_{j_-(\a, \c)}  (D)}  (\a^*U_{j_-(\a, \c)}  -m\frac{\chi}  {k}   )=
  \prod_{m=1}  ^{k}   (0-m\frac{\chi}  {k}  ) \times   \\
&\prod_{m=1}  ^{U_{j_-(\a, \c)}  (d)+\c^*U_{j_-(\a, \c)}  (D)}  (\c^*U_{j_-(\a, \c)}  -m\frac{\chi}  {k}  )   
, \\
&\text{The terms that come before the `` $\times $" sign are associated with the}  \\&\text{deformations along $T\C P^1_{\eps, \eps'} $.}  \\
\end{align*}  \begin{align*}   
&\text{For}   \ j\in \{j_+(\a, \b)   |  \a\sim\b\}    , \neq j_+(\a, \c)   , 
\\&  \ \ \prod_{m=1}  ^{U_j(d)+\a^*U_j(D)+{\tilde P}  ^{(\a, j_+(\a, \c)  ) }  (D)+{\tilde d}}  (\a^*U_j+{\tilde P}  ^{(\a, j_+(\a, \c)  ) }  -m\frac{\chi}  {k}  )   = \\
&  \prod_{m=1}  ^{k}  
 (\a^*U_j+{\tilde P}  ^{(\a, j_+(\a, \c)  ) }  -m\frac{\chi}  {k}  )   \times \\
&\prod_{m=1}  ^{U_j(d)+\c^*U_j(D)+{\tilde P}  ^{(\a, j_+(\a, \c)  ) }  (D)   -0+\a^*U_j(D)   -\c^*U_j(D)   
+{\tilde d}  -k }  (\c^*U_j-m\frac{\chi}  {k}  ), \\
&\text{for}   \ j\in \a\cap\c 
, \ \ \prod_{m=1}  ^{U_j(d)+\a^*U_j(D)}  (\a^*U_j-m\frac{\chi}  {k}  )   = \\
&  1 \times   
\prod_{m=1}  ^{U_j(d)+\c^*U_j(D)}  (\c^*U_j-m\frac{\chi}  {k}  ), \end{align*}  \begin{align*}  
&
\prod_{m=1}  ^{U_{j_+}  (d)+\a^*U_{j_+}  (D)+{\tilde P}  ^{(\a, j_+(\a, \c)  ) }  (D)+{\tilde d}}  (\a^*U_{j_+}  +{\tilde P}  ^{(\a, j_+(\a, \c)  ) }  -m\frac{\chi}  {k}  )   = \\
&  \prod_{m=1}  ^{k}  
 (0 -m\frac{\chi}  {k}  )   \times \\
&\prod_{m=1}  ^{U_{j_+}  (d)+\c^*U_{j_+}  (D)+{\tilde P}  ^{(\a, j_+(\a, \c)  ) }  (D)   -0+\a^*U_{j_+}  (D)+\c^*U_{j_+}  (D)   
+{\tilde d}  -k }  (\c^*U_{j_+}  -m\frac{\chi}  {k}  )   , \\
& \prod_{m=1}  ^{c_1(L)   (D)+{\tilde P}  ^{(\a, j_+(\a, \c)  ) }  (D)+{\tilde d}}  (c_1(L)+{\tilde P}  ^{(\a, j_+(\a, \c)  ) }  +m(-\chi/k)  )   \\ 
&= \prod_{m=1}  ^{k}  (c_1(L)+{\tilde P}  ^{(\a, j_+(\a, \c)  ) }  -m\frac{\chi}  {k}  )   \times \\
&\prod_{m=1}  ^{c_1(L)   (D)+\c^*U_j(D)+{\tilde P}  ^{(\a, j_+(\a, \c)  ) }  (D)   -0+\a^*U_j(D)   -\c^*U_j(D)+{\tilde d}  -k }  (c_1(L)+0-m\frac{\chi}  {k}  ), \\
&  (Q^Dq^{P^{\a}  (D) }  {\tilde q}  ^{{\tilde P}  ^{(\a, j_+(\a, \c)  ) }  (D) }  )    q^d {\tilde q}  ^{\tilde d}  = {\tilde q}  ^{-k}  q^{kd_{\a, \c}}  \times (Q^D q^{P^{\c}  (D) }  {\tilde q}  ^{0}  )   \\&q^{d-kd_{\a, \c}  +P^{\a}  (D)-P^{\c}  (D)}  {\tilde q}  ^{{\tilde d}  +k+{\tilde P}  ^{(\a, j_+(\a, \c)  ) }  (D)   -{\tilde P}  ^{\c}  (D) }  , \\
&\exp   \left( -\frac{(P^{\a}  t+{\tilde P}  ^{(\a, j_+(\a, \c)  ) }  {\tilde t}  )   k}  {\chi}  \right)   
\exp (dt+{\tilde d}  {\tilde t}  ) \\&   \exp (P^{\a}  (D)   t+{\tilde P}  ^{(\a, j_+(\a, \c)  ) }  (D)   {\tilde t}  )   =  1 \ \times \
\exp \left(-\frac{(P^{\c}  t+0{\tilde t}  )   k}  {\chi}  \right)  \\& 
\exp \left[dt+\left({\tilde d}  +k+{\tilde P}  ^{(\a, j_+(\a, \c)  ) }  (D)   -0+\a^*U_j(D)   -\c^*U_j(D)   \right)   {\tilde t}  \right]
\\&\exp (P^{\c}  (D)   t+0{\tilde t}  )   .
\end{align*}  In the last equality we use $({\tilde P}  ^{(\a, j_+(\a, \c)  ) }  -0)   /\chi=-1 $
from 2.2.

Factors on the R.H.S. which come before the multiplication
sign `` $\times $'' form the recursion coefficients
 ${\tilde q}  ^{-k}  Coeff_{(\a, j_+(\a, \c)  )   , \c}  (k)  $. Factors which
come after the multiplication sign form the term of the series $\J^{\c} $
evaluated at $z=-\chi/k $ and with the summation index $\tilde d $ replaced
with ${\tilde d}  +k+{\tilde P}  ^{(\a, j_+(\a, \c)  ) }  (D)   -0 $ and $d $ replaced with $d-kd_{\a, \b}  +\a^*U_j(D)   -\c^*U_j(D)  $. Reversing this change in the 
summation indices, we conclude that 
\[ \pi^{\perp}\Res_{z=\frac{(\a, j_+)   ^*{\tilde P}}  {k}}  \ \J^{(\a, j_+) }  (z)   \ dkz = i_A^* {\tilde q}  ^{-k}  q^{kd_{\a, \b}}  
Coeff_{(\a, j_+(\a, \c)  )   , \c}  (k)   \\\J^{\c}  (\frac{(\a, j_+)   ^*{\tilde P}}  {k}  )   , \]
as required.  Here we have used
\begin{align*}  \pi^{\perp}&\sum_{ {\tilde d}  \in \Z, 
D\in MC(B) }  \frac{({\tilde q} e^{\tilde t}  )   ^{\tilde d}  Q^D (i_A^*(-\frac{\chi}  {k})^{ht_{A  }  (D')} \J_{e(\cdot), L}  ^D(-\frac{\chi}{k}  , \tau) +\G^D_A  (-\frac{\chi}  {k}  , \tau))}  
{i_A^* \prod_{m=1}  ^{c_1(L)   (D)+{\tilde d}}  (c_1(L)   -m\frac{\chi}  {k}  ) }  =\\i_A^*&\sum_{ {\tilde d}  \in \Z, 
D\in MC(B) }  \frac{ \prod_{m=1}  ^{c_1(L)   (D) }  (c_1(L)   -m\frac{\chi}{k}  ) }  { \prod_{m=1}  ^{c_1(L)   (D) }  (c_1(L)   -m\frac{\chi}{k}  ) }  \times   \frac{({\tilde q} e^{\tilde t}  )   ^{\tilde d}  Q^D(-\frac{\chi}  {k}  )^{ht_{A  }  (D') } \J^D_{B}  (-\frac{\chi}  {k}  , \tau) }  
{ \prod_{m=1}  ^{\tilde d}  (c_1(L)   -c_1(L)   (D)\frac{\chi}  {k}     -m\frac{\chi}{k}  ) }.
\end{align*}

 \section{Mirrors}  
Our goal is to verify condition (1.b), (1.a)    of Theorem 2.

{\tt Proposition.}    For each element ${\tilde   \D}  \in  Ker \pi_* $, let $\Phi_{\tilde  \D}  (\tau _1, \dots, \tau_{B}  )  $ be a polynomial in coordinates $\tau_1, \dots, \tau_{B} $ along $H^2(B) $, with coefficients in $\Q(\l)  $.  Set $\Phi(\tau_1, \dots, \tau_{B}  ) :=\sum_{\tilde  \D\in Ker \pi_*}  (qe^t)^d({\tilde q} e^{\tilde t}  )^{\tilde d}  \Phi_{\tilde  \D}  (\tau_1, \dots, \tau_{B}  )  . $  Then, for each smooth family $\J(-z, \tau)   \subset  \LL_B $ ( $\subset   \LL_ A $ respectively)
\[e^{\Phi(z\p_{\tau_1}  , \dots, z\p_{\tau_{B}}  )/-z}  \J(-z, \tau)\]is contained in $\LL_B $ ( $ \LL_ A $ respectively). 

{\tt Proposition.}    For each element ${\tilde  \D}  \in  Ker \pi_* $, let   $\Phi_{\tilde  \D}  (x;\tau_1, \dots, \tau_{B}  )  $ be a polynomial in coordinates $\tau_1, \dots, \tau_{B} $ along $H^2(B) $ whose coefficients are smooth functions of $x $ valued in $\Q(\l)  $.  Set 
 $ \Phi(x;\tau_1, \dots, \tau_{B}  )  :=\sum_{\tilde \D\in Ker \pi_*}  (qe^t)^d \linebreak({\tilde q} e^{\tilde t}  )^{\tilde d}    \Phi_{\tilde  \D}  (x;\tau_1, \dots, \tau_{B}  )   . $ Suppose that, in a neighborhood of a given critical point of $\Phi(x)  $, the Taylor series of $\Phi(x)  $ converges, when $q, \tilde q, t, \tilde t, \l $ are valued in certain open sets of complex numbers.  Then, for each smooth family $\J(-z, \tau)   \subset  \LL_B $ ( $\subset   \LL_ A $ respectively)    the stationary phase asymptotics of the integral
\[ \int e^{\Phi(x;z\p_{\tau_1}  , \dots, z\p_{\tau_{B}}  )   /-z}  \J(-z, \tau)   dx\]is contained in $\LL_B $ ( $ \LL_ A $ respectively).

{\bf  7.2. The Quantum Riemann--Roch theorem.}  
Let $\LL $ and $\LL^{tw} $ be the overruled Lagrangian cones respectively:
of genus 0 Gromov--Witten theory of a target manifold $M $, and of such
a theory {\em twisted}   (in the sense of Example 1.5)    by a line bundle over $M $
with the equivariant 1st Chern class $\nu $. The cone $\LL $ lies in the 
symplectic loop space $(\H, \Omega)  $ based on the Poincar\'e pairing $(a, b)   =\int_M ab $, while $\LL^{tw} $ lies in $(\H, \Omega^{tw}  )  $ based on
 $(a, b)   ^{tw}  =\int_M ab/\nu $. The linear map $(\H, \Omega^{tw}  )    \to (\H, \Omega)  $ defined by $\f \mapsto \f/\sqrt{\nu} $ is a symplectomorphism.
The 
well-known asymptotics of the logarithm of the function

\[\Gamma(z, \nu)=\sqrt{\frac{\nu}  {2 \pi z}}  \int_0^{\infty} e^{(-x+\nu\ln x)   /z}  d\ln x=\sqrt{\frac{ \nu}  {2\pi z}}  {z^{\nu/z}}  \Gamma(\nu/z)   , \]
where $\Gamma(\nu/z)  $ is the Euler Gamma-function, is given by
\[ \widehat{\Gamma}  (z, \nu)    =
\exp\left\{\frac{-\nu+\nu \ln \nu}  {z}  +\sum_{m=1}  ^{\infty}  
\frac{B_{2m}}  {2m (2m-1) }  \left(\frac{z}  {\nu}  \right)   ^{2m-1}  \right\}    , \]where $B_{2m} $ are bernoulli numbers. 

{\tt Theorem}  (\cite{Coates-Givental}  ).\[ \LL = \sqrt{\nu^{-1}}  \widehat{\Gamma}  (z, \nu)   \LL^{tw}  .\]

In order to reach our goal, it therefore suffices to represent each of \begin{align*}    \left(\prod_{j\notin\c}  \sqrt{\c^*U_j^{-1}}  \widehat{\Gamma}  (-z, \c^*U_j)   \right) \c^*I_{E_{\a(A) }}  (z)&, \\ \left(\sqrt{\a^*U_{j_+}  ^{-1}}  \widehat{\Gamma}  (-z, \a^*U_{j_+}  )   \prod_{jj \neq j_+}  \sqrt{(\a, j_+)   ^*U_{A, jj}  ^{-1}}  \widehat{\Gamma}  (-z, (\a, j_+)   ^*U_{A, jj}  )   \right) \times &\\  (\a, {j_+}  )   ^*I_{E_{\a(A) }}  (z), \\ \text{and}   \left(\prod_{jj\neq [1, {\vec 0}  ]}  \sqrt{(\a, [1, {\vec 0}  ])   ^*U_{A, jj}  ^{-1}}  \widehat{\Gamma}  (-z, (\a, [1, {\vec 0}  ])   ^*U_{A, jj}  )   \right)   (\a, [1, {\vec 0}  ])  ^*I_{E_{\a(A) }}  (z)& \end{align*}  in the form described by the $-z\to z $ version of the Proposition.  
More precisely, we represent these series as \[ q^{-P^{\eps}  /z}   {\tilde q}  ^{-{\tilde P}  ^{\eps}  /z}  
\times (\text{the form in the $-z\to z $ version of the Proposition}  ).\]Given any $\rho\in H^2(B) $ and any scalar $\theta $, the multiplication by $e^{- \theta\rho /z} $ is theexponential of the operator $-\theta\p_{\rho}  +\sum_iQ_i \rho_i\p_{Q_i} $ which lies tangent to the Lagrangian cone of the base, by the Divisor equations.

{\bf 7.3.  Mirrors.}   Firstly, consider when the base is the point, in which case there is no $T $-fixed stratum $(\a, [1, {\vec 0}  ])  $. The generalization to arbitrary base  will follow from the  

{\tt Lemma.}   For each $\l\in  H^2(BT, \Q)  $, $\rho\in H^2(B, \Q)  $ and $\J $ satisfying the string and divisor equations, 
\[z\p_{\l+\rho}  \J(z, t)   =\sum_D(\l+\rho+z\rho(D)  )   Q^D\J^D(z, t)   .\]

{\em Proof.}   From the point of view of the cohomology of the base, $\l $ is a scalar. 
{\bf The $l $-blowup case.}  

Define $f:(\C^*)   ^{N+l}  \rightarrow \C $ to be the multivalued function
\[ f(x_1, \dots, x_N, y_1, \dots, y_l)   =\sum_{a=1}  ^ly_a+\sum_{j=1}  ^Nx_j+\l_j\ln x_j .\]Introduce
the complex submanifold 
\[ V= \{(x_1, \dots, x_N, y_1, \dots , y_l)   \ | \prod_{j=1}  ^N x_j^{{m}  _{ij}}  \prod_{a=1}  ^ly_a^{{m}  _{i(N+a) }}  =v_i, i=1, \dots, K+l 
\}     \]
of $(\C^*)   ^{N+l} $, where\[ ({m}  _{ij}  )   =
\left( \begin{array}  {ccc|ccc|ccc}  &&&&& &\\
   &(m^{\a}  _{ia}  )   &&& (m^{\check{\a}}  _{ib}  )   & &&0_{K\times l}  &\\
 &&&&&& &\\\hline
 0& \cdots& 0& 1&\cdots&1&&\\ \vdots&&\vdots&\vdots&&\vdots &&-\mathbb{I}  _l& \\0&\dots &0&1&\dots  &1&\\ \end{array}  \right)   , \]
 $m^{\a}  _{ia}  :=(m_{\a}  )   _{ij_a} $ for all $a=1, \dots, K $, and
\[v_i= \left\{ 
\begin{array}  {cll}  q_ie^{t_i}  & \text{ $i=1, \dots, K $}  \\
{\tilde q}  _{i-K} e^{{\tilde t}  _{i-K}}  & \text{ $i=K+1, \dots , K+l $}  .\end{array}  \right.       \]

When the base is the point and $l=1 $, these are the defining equations of the  blowup of the fiber at the fixed-point $\a $.
Thus, Corollary can be thought of as a bridge between the toric bundles theorem \cite{Brown}   and our Main Theorem.

In the sequel, consider the 1-blowup case  for ease of readability.

For each $\c\in E_{\a(A) }  ^T $, in connection with the Proposition, consider the oscillating integral

\begin{align*}   q^{-P^{\c}  /z}  \left(\frac{1}  {\sqrt{2 \pi z}}  \right)^{N-K}  &\times \\\int_{U_{\c}  \subset     V}  &e^{f(x, y;z\p_{\L}  )   / z}  
\frac{d\ln y}  {d\ln {\tilde q} e^{\tilde t}}  \frac{d\ln x_1 \w\cdots \w
d\ln x_N}  {d\ln q_1e^{t_1}  \w\cdots \w d\ln q_Ke^{t_K}}  \J_B(z, \tau)   , \end{align*}  
 where $U_{\c}  \subset  V $ is the non-compact cycle $\R^{N+1-(K+1) } $ parametrized by $\{x_j\}    _{j \notin\c} $. The differential operators may be processed at each order of  the series $\J(z)=\sum_DQ^D\J^D(z) $.  Let $x(\c)   \in U_{\c} $ be the critical point 
of $f|_{U_{\c}}  (x, y;\L) $ defined by the condition that its truncation modulo Novikov's variables is given by $x_j=\c^*U_j\ \forall j\notin\c $.
More precisely, for each $j\notin\c $, $x_j(\c)  $ admits a series expansion of the form $\c^* U_j+ {\mathcal O}  (\{q^d{\tilde q}  ^{\tilde d}   \}    _{{\tilde \D\in Ker \pi_*}   }  )  $ which solves the critical point equations, with the caveat
that in the case $\c={\tilde \a}  (B\setminus A) $ an extension of the Novikov ring of $E_{\a(A)} $ is needed for the expansion.

Multiply on both sides of the linear system\footnote{Here and in the sequel, the index range is implicit.  In particular, $i\leq K $ and $ b \leq N $ on the LHS. The system variables $\ln x_1, \dots, \ln x_N, \ln y $ are in correspondence with the columns of the LHS.}     \[\left( \begin{array}  {ccccc|c}   &&&&& 0\\
    && (m_{i b}  )   && &\vdots\\
   &&&&& 0 \\\hline
   0& \cdots& 0\ 1&\cdots&1&-1 \\ \end{array}  \right)   =\left( \begin{array}  {c}  \ln q_1e^{t_1}  \\
  \vdots\\
 \ln q_Ke^{t_K}  \\
\ln {\tilde q} e^{\tilde t}  \\ \end{array}  \right)   \]
by
 
\[\left( \begin{array}  {ccc|c}   &&& 0\\
   &(m^{\c}  _{ia}  )   ^{-1}  &&\vdots\\
 &&&0\\\hline
   0& \cdots& 0& 1\\ \end{array}  \right)   .\]
 This gives the linear system 
\[\left( \begin{array}  {cccccc|c}   &&&&&& 0\\
    &&&(m^{\c}  _{ia}  )   ^{-1}  (m_{i b}  )   &&& \vdots\\
   &&&&&& 0\\\hline
   0&\ \cdots&&0\ 1&\cdots&1&-1 \\ \end{array}  \right) =\left( \begin{array}  {c}  \sum_{i=1}  ^K(m_{\c}  ^{-1}  )   _{j_1i}  \ln q_ie^{t_i}  \\
     \vdots\\
  \sum_{i=1}  ^K(m_{\c}  ^{-1}  )   _{j_Ki}  \ln q_ie^{t_i}  \\
    \ln {\tilde q} e^{\tilde t}  \\ \end{array}  \right)   , \]whose $K\times K $ submatrix with columns $j_1, \dots, j_K\in\c $ in the upper $K $ rows is the identity matrix.

For each $\c\in E_{\a(A) }  ^T, j_a\in\c $ we have

\begin{align*}  \ln x_{j_a}  =-\sum_{j\notin\c}  \sum_{i=1}  ^K\left[(m_{\c}  ^{-1}  )   _{j_ai}  m_{ij}  \right]\ln  x_j
+&\\
\sum_{i=1}  ^K(m_{\c}  ^{-1}  )   _{j_ai}  \ln (q_ie^{t_i}  )   , & \end{align*}  \[ x_{j_a}  = \left(\prod_{i=1}  ^K (q_ie^{t_i}  )   ^{(m_{\c}  ^{-1}  )   _{j_ai}}  \right)   \prod_{j\notin\c}  x_j^{-\sum_{i=1}  ^K\left[(m_{\c}  ^{-1}  )   _{j_ai}  m_{i j}  \right]}  , \text {and}    \]

\[exp(x_{j_a}  /z)   =\sum_{n_{j_a}  =0}  ^ {\infty}  \frac{\left(\prod_{i=1}  ^K (q_ie^{t_i}  )   ^{n_{j_a}  (m_{\c}  ^{-1}  )   _{j_ai}}  \right)   \prod_{j\notin\c}  x_j^{-\sum_{i=1}  ^K\left[n_{j_a}  (m_{\c}  ^{-1}  )   _{j_ai}  m_{i j}  \right]}}  {n_{j_a}  !z^{n_{j_a}}}  .\]

Identify $n_{j_a}  \to \left\{ \begin{array}  {cl}  U_{j_a}  (d)\geq 0 & , j_a\in\c\cap\a\\U_{j_a}  (d)+{\tilde d}  \geq 0 &, j_a
\in\c\cap\a^c\end{array}  \right.        $ and identify $n\to -{\tilde d}  \geq 0 $. The RHS of the identification is nonnegative in view of the characterization  of the Mori cone of $E_{\a(A) } $ in section 6. 
The solution set to the above system of inequalities
defines a proper in general subset   ${\cal S}  \subset Ker\pi_* $ provided that $\c\neq \tilde \a(B\setminus A) $.

 Finally, do the elementary row operations that put a zero in the last row of the $j $th columns for each index value of $j =1, \dots, N $ for which $j\in\c $, 
and rewrite the exponential of \begin{align*}  y=({\tilde q} e^{\tilde t}  )   ^{-1}  \prod_{j\notin\a }  x_j\end{align*}   in terms of $\{x_j\}    _{j\notin\c} $:

\begin{align*}   &\exp(y/z)   =\sum_{n=0}  ^ {\infty}  \frac{y^n}  {n!z^n}  =\sum_{n=0}  ^ {\infty}  \frac{({\tilde q} e^{\tilde t}  ) ^{-n}  \prod_{i=1}  ^K(q_ie^{t_i}  )^{n\sum_{j_a\in\a^c\cap\c}  (m_{ia}  ^{\c}  )^{-1}}}  {n!z^n}  \times   \\&\left(\prod_{j\in\a\cap\c^c}  x_j^{-\sum_{j_a\in\c\cap\a^c}  \sum_{i=1}  ^K(m_{ia}  ^{\c}  )^{-1}  m_{ij}}  \prod_{j\notin\a \cup\c}  x_j^{1-\sum_{j_a\in\c\cap\a^c}  \sum_{i=1}  ^K(m_{ia}  ^{\c}  )^{-1}  m_{ij}}  \right)   ^n.\end{align*}

{Hence the factor $exp[(y+\sum_{j_a\in\c}  x_{j_a}  )   /z]$
of the integrand $e^{f/z} $ assumes the form \begin{align*} exp&(y/z)   \prod_{a=1}  ^Kexp(x_{j_a}  /z)   =\\ \sum_{(d , {\tilde d}  )\in   \cal S}  &\frac{(qe^t)   ^d({\tilde q} e^{\tilde t}  )   ^{\tilde d}}  {\prod_{m=1}  ^{-\tilde d}  mz\prod_{j\in\c\cap\a}  \prod_{m=1}  ^{U_j(d) }  mz
\prod_{j\in\a^c\cap\c}  \prod_{m=1}  ^{U_j(d)+\tilde d}  mz
}  \prod_{j\notin\c}  x_j^{-U_j(d) }  \prod_{j\notin\a \cup\c}  x_j^{-{\tilde d}}  .\end{align*}  Multiplying by $\l_{j_a} $ on the $\ln x_{j_a} $ formula and summing over $a=1\dots, K $ gives \begin{align*}  \sum_{a=1}  ^K\l_{j_a}  \ln (x_{j_a}  )   =-\sum_{j\notin\c}  \sum_{i, a=1}  ^K (m_{\c}  ^{-1}  )   _{j_ai}  \l_{j_a}  m_{ij}  \ln x_j+\sum _{i=1}  ^KP^{\c}  _i\ln (q_ie^{t_i}  )   =\\-\sum_{j\notin\c}  (\c^*U_j+\l_j)    \ln x_j+\sum _{i=1}  ^KP^{\c}  _i\ln (q_ie^{t_i}  ).\end{align*}    
Combining the ingredients of \[f|_{U_{\c}}  (x, y;\l)   =y+\sum_{j_a\in\c}  ^K(x_{j_a}  +\l_{j_a}  \ln x_{j_a}  )+\sum_{j\notin\c}  (x_j+\l_j \ln x_j), \]
 we obtain 
\begin{align*} e^{f|_{U_{\c}}  (x, y;\l)   /z}  =\\\prod_{i=1}  ^K(q_ie^{t_i}  )   ^{P_i^{\c}  /z}  \sum_{(d , {\tilde d}  )\in   \cal S}  \frac{(qe^t)   ^d({\tilde q} e^{\tilde t}  )   ^{\tilde d}}  {\prod_{m=1}  ^{-\tilde d}  mz\prod_{j\in\c\cap\a}  \prod_{m=1}  ^{U_j(d) }  mz\prod_{j\in\a^c\cap\c}  \prod_{m=1}  ^{U_j(d)+\tilde d}  mz}  \times \\ \prod_{j\notin\c} e^{x_j/z}  x_j^{-\c^*U_j/z-U_j(d)}  \prod_{j\notin\a \cup\c}  x_j^{-\tilde d}  .\end{align*}  By applying the Lemma, we deduce the differential operator version

\begin{align*} e^{f|_{U_{\c}}  (x, y;z\p_{\L}  )   / z}  \J_B(z, \tau)   =\sum_{D\in MC(B)}  \prod_{i=1}  ^K(q_ie^{t_i}  )   ^{P_i^{\c}  /z+P_i^{\c}  (D) }  \times &\\ \sum_{(d, {\tilde d}  )\in   \cal S}   \frac{Q^D(qe^t)   ^d({\tilde q} e^{\tilde t}  )   ^{\tilde d}}  {\prod_{m=1}  ^{-\tilde d}  mz\prod_{j\in\c\cap\a}  \prod_{m=1}  ^{U_j(d) }  mz\prod_{j\in\a^c\cap\c}  \prod_{m=1}  ^{U_j(d)+\tilde d}  mz}  \times \\\prod_{j\notin\c} e^{x_j/z}  x_j^{-\c^*U_j/z-\c^*U_j(D)   -U_j(d) }  &\prod_{j\notin\a \cup\c}  x_j^{-{\tilde d}}  \J_B^D(z, \tau)   .\end{align*}  Multiply the latter by $\prod_{j\notin\c}  d\ln x_j $ from the integrand, integrate by parts using the Gamma-function identities, and take stationary phase asymptotics to obtain
 \begin{align*}  \left(\prod_{j\notin\c}  \sqrt{\c^*U_j^{-1}}  \widehat{\Gamma}  (-z, \c^*U_j)   \right)    \c^*I^{pre}  _{E_{\a(A) }}  (z), \end{align*}   where
\begin{align*}  \c^*I^{pre}  _{E_{\a(A) }}  =e^{P^{\c}  t/z}   \sum_{D\in MC(B) }  \sum_{d\in \Z^K, {\tilde d}  \in\Z}  \\
\frac{(Q^Dq^{P^{\c}  (D) }  )    q^d{\tilde q}  ^{{\tilde d}} e^{{\tilde d}  {\tilde t}} e^{dt} e^{P^{\c}  (D)   t}}  
{ \prod_{j\in \{j_+(\a, \b)   |  \a\sim\b\}  }  \prod_{m=1}  ^{U_j(d)+\c^*U_j(D)+{\tilde d}}  (\c^*U_j+mz) }  \times \\
\frac{1}  {\prod_{j\notin \{j_+(\a, \b)   |  \a\sim\b\}  }  \prod_{m=1}  ^{U_j(d)+\c^*U_j(D) }  (\c^*U_j+mz) }  \times \\ \frac{1}  {\prod_{m=1}  ^{-{\tilde d}}  (mz)}  \J_B^D(z, \tau). \end{align*}  
The series $\c^*I_{E_{\a(A)}} $ of the Main Theorem is related to $\c^*I_{E_{\a(A)}}  ^{pre} $ by 
\begin{align*}  \c^*I_{E_{\a(A) }}  ^{D, {\tilde d}  , d}  (z, t, {\tilde t}  , \tau, q, {\tilde q}  , Q)   &=\\z^{ht_A(D')}  \times  &\frac{(\c^*I^{pre}  _{E_{\a(A) }}  )^{D, {\tilde d}  , d}  (z, t, {\tilde t}  , \tau, q, {\tilde q}  , Q) }  {\prod_{m=1}  ^{{\tilde d}}  (c_1(L)+c_1(L)   (D)   z+mz) }  , \end{align*} $\ \forall D \leq D' $.

{\tt Lemma.}      Given any complex line bundle $L \to B $, the series 

\begin{align*}  \sum_{ {\tilde d}  \in \Z, 
D\in MC(B) }  \frac{({\tilde q} e^{\tilde t}  )   ^{\tilde d}  Q^D \J^D_{B}  (z, \tau) }  
{ \prod_{m=1}  ^{\tilde d}  (c_1(L)+c_1(L)   (D)   z+mz) }  
\end{align*}  
is contained in $\LL_B $ for all  formal values of ${\tilde t}  , \tau, {\tilde q}  , Q $. 

{\em Proof.}    The cone $\LL_B $ is preserved by both numerator and denominator in the ratio of operators on the LHS of
\begin{align*}  
\sum_{ {\tilde d}  \in \Z}  ({\tilde q} e^{\tilde t}  )   ^{\tilde d}  
\frac{\left(\int e^{x/z}  x^{\frac{z\p_{c_1(L) }}  {z}}  d\ln x\right)   ^{\widehat{}}}  {\left(\int e^{x/z}  x^{\frac{z\p_{c_1(L) }}  {z}  +{\tilde d}}  d\ln x\right)   ^{\widehat{}}}  &\J_{B}  (z, \tau)   
=\\\sum_{ {\tilde d}  \in \Z, 
D\in MC(B) }  &\frac{({\tilde q} e^{\tilde t}  )   ^{\tilde d}  Q^D \J^D_{B}  (z, \tau) }  
{ \prod_{m=c_1(L)   (D)+1}  ^{c_1(L)   (D)+{\tilde d}}  (c_1(L)+mz) }  =\\
\sum_{ {\tilde d}  \in \Z, 
D\in MC(B) }  &\frac{({\tilde q} e^{\tilde t}  )   ^{\tilde d}  Q^D \J^D_{B}  (z, \tau) }  
{ \prod_{m=1}  ^{\tilde d}  (c_1(L)+c_1(L)   (D)   z+mz) }  , \end{align*}   
where $\widehat{} $ denotes stationary phase asymptotics.

Henceforth denote $\J_A^{D, D'}(z, \tau)={\bf z^{ht_{A  }  (D') }} i_A^*\J^D_{e(\cdot), L}  (z, \tau) +\G^D_{ A}  (z, \tau) $.
For each $(\a, j_+)   \in E_{\a(A) }  ^T $ in connection with the Proposition, consider
the oscillating integral
\begin{align*}  {\tilde q}  ^{-{\tilde P}  ^{(\a, j_+) }  /z}  q^{-P^{\a}  /z}  &\left(\frac{1}  {\sqrt{2 \pi z}}  \right)^{N-K+1}  \int_Uexp\left(\frac{x'-\ln x'(z\p_{c_1(L)}  +z\p_{\tilde t}  )}  {z}  \right)d\ln x'\times \\ \int_{U_{\a, j_+}  \times U\subset     V}  &e^{f(x, y;z\p_{\L}  )   / z}  \frac{d\ln y}  {d\ln {\tilde q} e^{\tilde t}}   \frac{d\ln x_1 \w\cdots \w
d\ln x_N}  {d\ln q_1e^{t_1}  \w\cdots \w d\ln q_Ke^{t_K}}  \J_{A}  (z, \tau)   , \end{align*}   where $U_{\a, j_+}   \times U\subset  V\times \C $ is  the non-compact cycle $\R^{N+1-(K+1) }  \times \R $ parametrized by $\{y\}    \cup\{x_j\}    _{j \notin\a , \neq j_+}  \cup\{x'\}   $.
 The differential operators may be processed at each order of the series $\J(z)=\sum_DQ^D\J^D(z) $.  Then for each $a\in H^*(A) $ and $b\in H^*(B) $ interpret $ab $ to mean $a(\a, j_+)^*b $.  Let $x(\a, j_+)   \in U_{\a, j_+}  \times U $ be the critical point of \[x'-\ln x'(c_1(L)+{\tilde P}  ^{(\a, j_+)}  )+f|_{U_{\a, j_+}  \times U}  (x, y;\L)\]
defined by the condition that its truncation modulo Novikov's variables is given by $x_j=(\a, j_+)^*U_{A, j}  \ \forall j\notin\a , \neq j_+ $, by $y=(\a, j_+)^*U_{j_+} $, and by $x'=(\a, j_+)^*U_{A, [1, \vec 0]}    $.
More precisely, for each $j\notin\a , \neq j_+ $, $x_j(\a, j_+)  $ admits a series expansion of the form $(\a, j_+)^* U_{A, j}  + {\mathcal O}  (\{q^d{\tilde q}  ^{\tilde d}   \}    _{{\tilde \D\in Ker \pi_*}   }  )  $ which solves the critical point equations.  Similarly for $x' $ and $y $.

Multiply on both sides of the linear system \[\left( \begin{array}  {ccc|ccc|c}  &&&&&& 0\\
   & (m^{\a}  _{ia}  )   &&&(m^{\check{\a}}  _{ib}  )   & &\vdots\\
 & &&&&& 0 \\\hline
   0& \cdots& 0& 1&\cdots&1&-1 \\ \end{array}  \right)   = \left( \begin{array}  {c}  \ln q_1e^{t_1}  \\
    \vdots\\ 
   \ln q_Ke^{t_K}  \\
   \ln {\tilde q} e^{\tilde t}  \\ \end{array}  \right)   \]
by 
\[\left( \begin{array}  {ccc|c}   &&& 0\\
   &(m^{\a}  _{ia}  )   ^{-1}  &&\vdots\\
 &&&0\\\hline
   0& \cdots& 0& 1\\ \end{array}  \right)   .\]
This gives the linear system 
\[\left( \begin{array}  {ccc|ccc|c}  &&&&&& 0\\
   & 
\mathbb{I}  
&&&(m^{\a}  _{ia}  )   ^{-1}  (m^{\check{\a}}  _{ib}  )   & &\vdots\\
 & &&&&& 0\\\hline
   0& \cdots& 0& 1&\cdots&1&-1 \\ \end{array}  \right)    = \left( \begin{array}  {c}  \sum_{i=1}  ^K(m_{\a}  ^{-1}  )   _{j_1i}  \ln q_ie^{t_i}  \\
     \vdots\\
 \sum_{i=1}  ^K(m_{\a}  ^{-1}  )   _{j_Ki}  \ln q_ie^{t_i}  \\
   \ln {\tilde q} e^{\tilde t}  \\ \end{array}  \right)   .\]

Perform the $K $ elementary row operations on both sides of the equation which result in a $0 $ in the first $K $ rows of the $j_+ $ column of the LHS.  Then for each $j_a\in \a $, 

\begin{align*}  x_{j_a}  =\left(\prod_{i=1}  ^K (q_ie^{t_i}  )   ^{(m_{\a}  ^{-1}  )   _{j_ai}}  ({\tilde q} e^{\tilde t}  )   ^{-
(m_{\a}  ^{-1}  )   _{j_ai}  m_{ij_+}}  \right)   y^{-(m_{\a}  ^{-1}  )   _{j_ai}  m_{ij_+}}  \times \\ \prod_{j\notin\a }  x_j^{-\sum_{i=1}  ^K\left[(m_{\a}  ^{-1}  )   _{j_ai}  m_{i j}  -(m_{\a}  ^{-1}  )   _{j_ai}  m_{ij_+}  \right]}  , \text {and}   \end{align*}  
\begin{align*}  exp(x_{j_a}  /z)   =\sum_{n_{j_a}  =0}  ^ {\infty}  \frac{\left(\prod_{i=1}  ^K (q_ie^{t_i}  )   ^{n_{j_a}  (m_{\a}  ^{-1}  )   _{j_ai}}  ({\tilde q} e^{\tilde t}  )   ^{-
n_{j_a}  (m_{\a}  ^{-1}  )   _{j_ai}  m_{ij_+}}  \right) }  {n_{j_a}  !z^{n_{j_a}}}  \times \\ y^{-n_{j_a}  (m_{\a}  ^{-1}  )   _{j_a}  m_{ij_+}}  \prod_{j\notin\a }  x_j^{-\sum_{i=1}  ^Kn_{j_a}  \left[(m_{\a}  ^{-1}  )   _{j_ai}  m_{i j}  -(m_{\a}  ^{-1}  )   _{j_ai}  m_{ij_+}  \right]}  .\end{align*}  
Next, \begin{align*}   x_{j_+}  ={\tilde q} e^{\tilde t}  y\prod_{j\notin\a , \neq j_+}  x_j^{-1}  , \text{and}  \end{align*}  \[exp(x_{j_+}  /z)   =\sum_{n=0}  ^ {\infty}  \frac{x_{j_+}  ^n}  {n!z^n}  =\sum_{n=0}  ^ {\infty}  \frac{({\tilde q} e^{\tilde t}  y\prod_{j\notin\a , \neq j_+}  x_j^{-1}  )   ^n}  {n!z^n}  .\]

 Identify $n_{j_a}  \to U_{j_a}  (d)\geq 0 \forall a=1, \dots, K $ and identify $n\to U_{j_+}  (d)+{\tilde d}  \geq 0 $.  The RHS of the identification is nonnegative in view of the characterization  of the Mori cone of $E_{\a(A) } $ in section 6.  The solution set to the above system of inequalities
defines a proper in general subset   ${\cal S}  \subset Ker\pi_* $.{
\begin{align*} exp(x_{j_
+}  /z)   \prod_{a=1}  ^Kexp&(x_{j_a}  /z)   =\\ \sum_{(d, {\tilde d}  )\in   \cal S }  &\frac{(qe^t)   ^d({\tilde q} e^{\tilde t}  )   ^{\tilde d}}  {\prod_{m=1}  ^{U_{j_+}  (d)+\tilde d}  mz\prod_{j\in\a, \neq j_+}  \prod_{m=1}  ^{U_j(d) }  mz}  y^{U_{j_+}  (d)+\tilde d}  \prod_{j\notin\a }  x_j^{-U_j(d)   -{\tilde d}}  .\end{align*}  For each $ j_a\in\a $ we have

\begin{align*}  \ln x_{j_a}  =-\sum_{j\notin\a , \neq j_+}  \sum_{i=1}  ^K\left[(m_{\a}  ^{-1}  )   _{j_ai}  m_{ij}  -(m_{\a}  ^{-1}  )   _{j_ai}  m_{ij_+}  \right]\ln  x_j
+&\\
\sum_{i=1}  ^K(m_{\a}  ^{-1}  )   _{j_ai}  \ln (q_ie^{t_i}  )-\sum_{i=1}  ^K(m_{\a}  ^{-1}  )   _{j_ai}  m_{ij_+}  \ln  y-\sum_{i=1}  ^K(m_{\a}  ^{-1}  )   _{j_ai}  &m_{ij_+}  \ln  {\tilde q} e^{\tilde t}  .\end{align*}   
Thus, we arrive at
\begin{align*}  \sum_{j=1}  ^N(x_j&+\ln x_j)+ \\  \sum_{j\notin\a , \neq j_+}  \l_j&\ln x_j-\sum_{j\notin\a , \neq j_+}  \sum_{j_a\in\a}  \sum_{i=1}  ^K\l_{j_a}  \left[(m_{\a}  ^{-1}  )   _{j_ai}  m_{ij}  -(m_{\a}  ^{-1}  )   _{j_ai}  m_{ij_+}  \right]\ln x_j
+\\\sum_{j_a\in\a}  &\sum_{i=1}  ^K\l_{j_a}  (m_{\a}  ^{-1}  )   _{j_ai}  \ln (q_ie^{t_i}  )   -\sum_{j_a\in\a}  \sum_{i=1}  ^K\l_{j_a}  (m_{\a}  ^{-1}  )   _{j_ai}  m_{ij_+}  \ln  y-\\  \sum_{j_a\in\a}  &\sum_{i=1}  ^K\l_{j_a}  (m_{\a}  ^{-1}  )   _{j_ai}  m_{ij_+}  \ln  {\tilde q} e^{\tilde t}  .\end{align*}   
For each $j\notin\a , \neq j_+ $, 
\begin{align*}   \left(\l_j-\sum_{j_a\in\a}  \sum_{i=1}  ^K\l_{j_a}  \left[(m_{\a}  ^{-1}  )   _{j_ai}  m_{ij}  -(m_{\a}  ^{-1}  )   _{j_ai}  m_{ij_+}  \right]\right)   \ln x_j =&\\ -(\a, j_+)   ^*U_{A, j}  \ln x_j&, 
\end{align*}   while the third and fourth lines of the former is given by  \begin{align*}  \sum_{i=1}  ^KP^{\a}  _i \ln q_ie^{t_i}  - \a^*U_{j_+}  \ln y-\a^*U_{j_+}  \ln {\tilde q} e^{\tilde t}  .\end{align*}   

Combining  the ingredients of \begin{align*}  f|_{U_{\a, j_+}  \times U}  (x, y;\l)   =&y+\\(x_{j_+}  +\ln x_{j_+}  )+&\sum_{j_a\in\a, \neq j_+}  (x_{j_a}  +\l_{j_a}  \ln x_{j_a}  )+\sum_{j\notin\, a, \neq j_+}  (x_j+\l_j \ln x_j), \end{align*}  we obtain 
\begin{align*} e^{f|_{U_{\a, j_+}  \times U}  (x, y;\l)   /z}  =({\tilde q} e^{{\tilde t}}  )   ^{{\tilde P}  ^{(\a, j_+) }  /z}  \prod_{i=1}  ^K(q_ie^{t_i}  )   ^{P_i^{\a}  /z}  \times \\ \sum_{(d, {\tilde d}  )\in   \cal S }  \frac{(qe^t)   ^d({\tilde q} e^{\tilde t}  )   ^{\tilde d}}  {\prod_{m=1}  ^{U_{j_+}  (d)+\tilde d}  mz\prod_{j\in\a}  \prod_{m=1}  ^{U_j(d) }  mz}  \times \\e^{y/z}  y^{{\tilde P}  ^ {(\a, j_+)}  /z+\tilde d}  \prod_{j\notin\a , \neq j_+} e^{x_j/z}  x_j^{-\a^*U_j/z-{\tilde P}  ^{(\a, j_+)}  /z -U_j(d)   -{\tilde d}}  .\end{align*}  By applying the Lemma, we deduce the differential operator version   
\begin{align*} e^{f|_{U_{\a, j_+}  \times U}  (x, y;z\p_{\L}  )   / z}  \J_{A}  (z, \tau)    =\\\sum_{D\in MC(B)}  ({\tilde q} e^{{\tilde t}}  )   ^{{\tilde P}  ^{(\a, j_+) }  /z+{\tilde P}  ^{(\a, j_+) }  (D) }  \prod_{i=1}  ^K(q_ie^{t_i}  )   ^{P_i^{\a}  /z+P_i^{\a}  (D) }  \times \\ \sum_{(d, {\tilde d}  )\in    \cal S }  \frac{Q^D(qe^t)   ^d({\tilde q} e^{\tilde t}  )   ^{\tilde d}}  {\prod_{m=1}  ^{
U_{j_+}  (d)+\tilde d}  mz\prod_{j\in\a}  \prod_{m=1}  ^{U_j(d) }  mz} e^{y/z}  y^{{\tilde P}  ^ {(\a, j_+)}  /z+{\tilde P}  ^{(\a, j_+)}  (D)+\tilde d}  \times \\\prod_{j\notin\a , \neq j_+} e^{x_j/z}  x_j^{-\a^*U_j/z
-\a^*U_j(D)-U_j(d)-{\tilde P}  ^{(\a, j_+)}  /z  -{\tilde P}  ^{(\a, j_+) }  (D)   -{\tilde d}}  \J_A^{D, D'}(z, \tau)   .\end{align*}  
Apply the operator \[exp\left(\frac{x'-\ln x'(z\p_{c_1(L)}  +z\p_{\tilde t}  )}  {z}  \right)\]
 to the latter to obtain the additional factor \begin{align*}  x'^{-c_1(L)/z-{\tilde P}  ^{(\a, j_+)}  /z-c_1(L)(D)-{\tilde P}  ^{(\a, j_+)}  (D)-{\tilde d}}  .\end{align*}  
Multiply the latter by $d\ln x'\ d\ln y \prod_{j\notin\a , \neq j   _+}  d\ln x_j $ from the integrand, integrate by parts using the Gamma-function identities, and take stationary phase asymptotics to obtain

  \begin{align*}  \left(\sqrt{\a^*U_{j_+}  ^{-1}}  \widehat{\Gamma}  (-z, \a^*U_{j_+}  )   \prod_{jj \neq j_+}  \sqrt{(\a, j_+)   ^*U_{A, jj}  ^{-1}}  \widehat{\Gamma}  (-z, (\a, j_+)   ^*U_{A, jj}  )   \right)\times    \\   (\a, {j_+}  )   ^*I  _{E_{\a(A) }}  (z), \end{align*}   where
\begin{align*}  (\a, j_+(\a, \c)  )   ^*I  _{E_{\a(A) }}  =e^{{\tilde P}  ^{(\a, j_+(\a, \c)  ) }  {\tilde t}  /z} e^{P^{\a}  t/z}   
\sum_{D\in MC(B) }  
\sum_{d \in \Z^K, {\tilde d}  \in\Z}  
\\
 \frac{ (Q^Dq^{{\tilde P}  ^{(\a, j_+(\a, \c)  ) }  (D) }  q^{P^{\a}  (D) }  )   q^d{\tilde q}  ^{{\tilde d}} e^{{\tilde d}  {\tilde t}} e^{dt} e^{{\tilde P}  ^{(\a, j_+(\a, \c)  ) }  (D)   t} e^{P^{\a}  (D)   t}}  
{\prod_{j\notin \{j_+(\a, \b)   |  \a\sim\b\}  }  
\prod_{m=1}  ^{U_j(d)+\a^*U_j(D) }  (\a^*U_j+mz) }  \times \\
\frac{1}  
{ \prod_{j\in \{j_+(\a, \b)   |  \a\sim\b\}  }  
\prod_{m=1}  ^{U_j(d)+\a^*U_j(D)+{\tilde P}  ^{(\a, j_+(\a, \c)  ) }  (D)+{\tilde d}}  (\a^*U_j+{\tilde P}  ^{(\a, j_+(\a, \c)  ) }  +mz) }  \times \\
\frac{1}  {\prod_{m=1}  ^{-{\tilde d}  -{\tilde P}  ^{(\a, j_+(\a, \c)  ) }  (D) }  (-{\tilde P}  ^{(\a, j_+(\a, \c)  ) }  +mz)   
}  \times \\ \frac{1}  {\prod_{m=1}  ^{c_1(L)   (D)+{ \tilde P}  ^{(\a, j_+(\a, \c)  ) }  (D)+{\tilde d}}  (c_1(L)+{\tilde P}  ^{(\a, j_+(\a, \c)  ) }  +mz) }  
\J_A^{D, D'}  (z, \tau)   .\end{align*}  
 $\ \forall D \leq D' $.

Let $g(x, y;z\p_{\L}  )  $ denote the expression
\begin{align*}   y+y{\tilde q} e^{\tilde t}  \prod_{j\notin\a }  x_j^{-1}  +\sum_{j\notin\a }  x_j+\ln x_j(z\p_{\a^*U_j}  
+z\p_{{\tilde P}  ^{(\a, [1, \vec 0]) }}  )+\sum_{i= 1}  ^K\ln(q_ie^{t_i}  )z\p_{P_i^{\a}}  +\\(\ln ({\tilde  q} e^{\tilde t}  )-\ln y) 
z\p_{{\tilde P}  ^{(\a, [1, \vec 0])}}  +\prod_{j\notin\a }  \left(\prod_{i=1}  ^K (q_ie^{t_i}  )   ^{(m_{\a}  ^{-1}  )   _{j_ai}}  \right)x_j^{-\sum_{i=1}  ^K(m_{\a}  ^{-1}  )   _{j_ai}  m_{i j}}  .\end{align*}  
Let $(d, \tilde d) $ index the solution set of degrees in section 6, in the case $\eps=(\a, [1, \vec 0]) $.  Then, 
\begin{align*} e^{\frac{g(x, y;z\p_{\L}  ) }  {z}}  \J_A(z, \tau)    =\\\sum_{D\in MC(B)}  ({\tilde q} e^{{\tilde t}}  )   ^{{\tilde P}  ^{(\a, [1, \vec 0]) }  /z+{\tilde P}  ^{(\a, [1, \vec 0]) }  (D) }  \prod_{i=1}  ^K(q_ie^{t_i}  )   ^{P_i^{\a}  /z+P_i^{\a}  (D) }  \times \\ \sum_{(d, {\tilde d}  )}  \frac{Q^D(qe^t)   ^d({\tilde q} e^{\tilde t}  )   ^{\tilde d}}  {\prod_{m=1}  ^{\tilde d}  mz\prod_{j\in\a}  \prod_{m=1}  ^{U_j(d) }  mz} e^{y/z}  y^{{\tilde P}  ^{(\a, [1, \vec 0])}  /z+{\tilde P}  ^{(\a, [1, \vec 0])}  (D)+\tilde d}  \times \\\prod_{j\notin\a } e^{x_j/z}  x_j^{-\a^*U_j/z-\a^*U_j(D)-U_j(d)-{\tilde P}  ^{(\a, [1, \vec 0])}  /z      -{\tilde P}  ^{(\a, [1, \vec 0]) }  (D)   -{\tilde d}}  \J_A^{D, D'}(z, \tau)   .\end{align*}  

Multiply by $d \ln y\ \prod_{j\notin\a }  d\ln x_j $ from the integrand, integrate by parts using the Gamma-function identities, and take stationary phase asymptotics to obtain
\begin{align*}  \left(\prod_{jj\neq [1, {\vec 0}  ]}  \sqrt{(\a, [1, {\vec 0}  ])   ^*U_{A, jj}  ^{-1}}  \widehat{\Gamma}  (-z, (\a, [1, {\vec 0}  ])   ^*U_{A, jj}  )   \right)   (\a, [1, {\vec 0}  ])   ^*I _{E_{\a(A) }}  (z), \ \end{align*}   where
\begin{align*}  (\a, [1, {\vec 0}  ])   ^*I_{E_{\a(A) }}  =e^{{\tilde P}  ^{(\a, [1, {\vec 0}  ]) }  {\tilde t}  /z} e^{P^{\a}  t/z}  
\sum_{D\in MC(B) }  
\sum_{d\in \Z^K, {\tilde d}  \in\Z}  & \\
\frac{ (Q^D{\tilde q}  ^{{\tilde P}  ^{(\a, [1, {\vec 0}  ]) }  (D) }  q^{P^{\a}  (D) }  )    q^d{\tilde q}  ^{{\tilde d}} e^{{\tilde d}  {\tilde t}}  
e^{dt} e^{{\tilde P}  ^{(\a, [1, {\vec 0}  ]) }  (D)   t} e^{P^{\a}  (D)   t}}  
{\prod_{j\notin \{j_+(\a, \b)   |  \a\sim\b\}  }  
\prod_{m=1}  ^{U_j(d)+\a^*U_j(D) }  (\a^*U_j+mz) }  
& \times \\\frac{1}  { \prod_{j\in \{j_+(\a, \b)   |  \a\sim\b\}  }  
\prod_{m=1}  ^{U_j(d)+\a^*U_j(D)+{\tilde P}  ^{(\a, [1, {\vec 0}  ]) }  (D)+{\tilde d}}  \ (\a^*U_j+{\tilde P}  ^{(\a, [1, {\vec 0}  ]) }  +mz) }  & \times \\
 \frac{1}  {\prod_{m=1}  ^{\tilde d}  (mz)   \prod_{m= 1}  ^{c_1(L)   (D)   -{\tilde d}}  (c_1(L)+mz) }   \J_A^{D, D'}(z, \tau)   .\end{align*}  

{\bf Acknowledgements.}   This research was mostly completed while JB was a postdoc at IBS Center for Geometry and Physics.  JB thanks the entire Postech community.  JB is also thankful for invitation to  sketch the results here at the Workshop on Differential Geometry in Sept. 2014 at the IMS, Chinese University of Hong Kong.

{\em E-mail:}   {\tt brown.jeff08@gmail.com}  
\end{document}